

\def\LaTeX{%
  \let\Begin\begin
  \let\End\end
  \let\salta\relax
  \let\finqui\relax
  \let\futuro\relax}

\def\UK{\def\our{our}\let\sz s}
\def\USA{\def\our{or}\let\sz z}

\UK



\LaTeX

\USA


\salta

\documentclass[twoside,12pt]{article}
\setlength{\textheight}{24cm}
\setlength{\textwidth}{16cm}
\setlength{\oddsidemargin}{2mm}
\setlength{\evensidemargin}{2mm}
\setlength{\topmargin}{-15mm}
\parskip2mm


\usepackage[usenames,dvipsnames]{color}
\usepackage{amsmath}
\usepackage{amsthm}
\usepackage{amssymb}
\usepackage[mathcal]{euscript}

\usepackage{cite}
%
%


\definecolor{viola}{rgb}{0.3,0,0.7}
\definecolor{ciclamino}{rgb}{0.5,0,0.5}

\def\previousgianni #1{{\color{ciclamino}#1}}
\def\juerg #1{{\color{PineGreen}#1}}
\def\gianni #1{{\color{red}#1}}
\def\pier #1{{\color{red}#1}}
\def\newpier #1{{\color{blue}#1}}
\def\jupiter #1{{\color{red}#1}}

\def\revis #1{{\color{red}#1}}
\def\revis #1{#1}

\def\previousgianni #1{#1}
\def\juerg #1{#1}
\def\gianni #1{#1}
\def\pier #1{#1}
\def\newpier #1{#1}
\def\jupiter #1{#1}




\bibliographystyle{plain}


%

\finqui

\def\Beq{\Begin{equation}}
\def\Eeq{\End{equation}}
\def\Bsist{\Begin{eqnarray}}
\def\Esist{\End{eqnarray}}

\def\Bthm{\Begin{theorem}}
\def\Ethm{\End{theorem}}

\def\Bprop{\Begin{proposition}}
\def\Eprop{\End{proposition}}

\def\Brem{\Begin{remark}\rm}
\def\Erem{\End{remark}}

\def\Bnot{\Begin{notation}\rm}
\def\Enot{\End{notation}}
\def\Bdim{\Begin{proof}}
\def\Edim{\End{proof}}
\def\Bcenter{\Begin{center}}
\def\Ecenter{\End{center}}
\let\non\nonumber




\def\step #1 \par{\medskip\noindent{\bf #1.}\quad}


\def\Lip{Lip\-schitz}

\def\aand{\quad\hbox{and}\quad}

\def\lhs{left-hand side}
\def\rhs{right-hand side}
\def\sfw{straightforward}



\def\multibold #1{\def\arg{#1}%
  \ifx\arg\pto \let\next\relax
  \else
  \def\next{\expandafter
    \def\csname #1#1#1\endcsname{{\bf #1}}%
    \multibold}%
  \fi \next}

\def\pto{.}

\def\multical #1{\def\arg{#1}%
  \ifx\arg\pto \let\next\relax
  \else
  \def\next{\expandafter
    \def\csname cal#1\endcsname{{\cal #1}}%
    \multical}%
  \fi \next}


\def\multimathop #1 {\def\arg{#1}%
  \ifx\arg\pto \let\next\relax
  \else
  \def\next{\expandafter
    \def\csname #1\endcsname{\mathop{\rm #1}\nolimits}%
    \multimathop}%
  \fi \next}

\multibold
qwertyuiopasdfghjklzxcvbnmQWERTYUIOPASDFGHJKLZXCVBNM.

\multical
QWERTYUIOPASDFGHJKLZXCVBNM.

\multimathop
diag dist div dom mean meas sign supp .


\def\accorpa #1#2{\eqref{#1}--\eqref{#2}}
\def\Accorpa #1#2 #3 {\gdef #1{\eqref{#2}--\eqref{#3}}%
  \wlog{}\wlog{\string #1 -> #2 - #3}\wlog{}}


\def\separa{\noalign{\allowbreak}}

\def\somma #1#2#3{\sum_{#1=#2}^{#3}}

\def\<#1>{\mathopen\langle #1\mathclose\rangle}
\def\norma #1{\mathopen \| #1\mathclose \|}
\def\Norma #1{\Bigl\| #1 \Bigr\|}

\def\[#1]{\mathopen\langle\!\langle #1\mathclose\rangle\!\rangle}

\def\iot {\int_0^t}
\def\ioT {\int_0^T}

\def\intQ{\int_Q}
\def\iO{\int_\Omega}

\def\dt{\partial_t}
\def\dn{\partial_\nu}

\def\cpto{\,\cdot\,}

\def\checkmmode #1{\relax\ifmmode\hbox{#1}\else{#1}\fi}

\def\aeQ{\checkmmode{a.e.\ in~$Q$}}

\def\aet{\checkmmode{a.e.\ in~$(0,T)$}}

\def\aat{\checkmmode{for a.a.~$t\in(0,T)$}}


\def\erre{{\mathbb{R}}}

\def\enne{{\mathbb{N}}}




\def\genspazio #1#2#3#4#5{#1^{#2}(#5,#4;#3)}
\def\spazio #1#2#3{\genspazio {#1}{#2}{#3}T0}

\def\L {\spazio L}
\def\H {\spazio H}
\def\W {\spazio W}


\def\Lx #1{L^{#1}(\Omega)}
\def\Hx #1{H^{#1}(\Omega)}

\def\LQ #1{L^{#1}(Q)}

\def\Luno{\Lx 1}
\def\Ldue{\Lx 2}

\def\Hdue{\Hx 2}
\def\Hunoz{{H^1_0(\Omega)}}


\def\LQ #1{L^{#1}(Q)}


\let\theta\vartheta

\let\phi\varphi

\let\hat\widehat

\let\TeXchi\chi                         
\newbox\chibox
\setbox0 \hbox{\mathsurround0pt $\TeXchi$}
\setbox\chibox \hbox{\raise\dp0 \box 0 }
\def\chi{\copy\chibox}



\def\Az #1{A_0^{#1}}
\def\Vz #1{V_0^{#1}}
\def\VA #1{V_A^{#1}}
\def\VB #1{V_B^{#1}}

\def\Beta{\hat\beta}
\def\betal{\beta_\lambda}
\def\Betal{\hat\beta_\lambda}
\def\betaz{\beta^\circ}
\def\Pi{\hat\pi}
\def\Lpi{L_\pi}

\def\yl{y^\lambda}
\def\mul{\mu^\lambda}

\def\yln{y^{\lambda_n}}
\def\muln{\mu^{\lambda_n}}
\def\Betaln{\Beta_{\lambda_n}}

\def\Betalm{\Beta_{\lambda_m}}

\def\yh{\hat y_h}
\def\muh{\hat\mu_h}
\def\overyh{\overline y_h}
\def\overmuh{\overline\mu_h}
\def\underyh{\underline y_h}
\def\undermuh{\underline\mu_h}
\def\overuh{\overline u_h}

\def\yn{y^n}
\def\ynp{y^{n+1}}
\def\ynm{y^{n-1}}
\def\mun{\mu^n}
\def\munp{\mu^{n+1}}
\def\munm{\mu^{n-1}}
\def\un{u^n}
\def\unp{u^{n+1}}
\def\dhyn{\frac{\ynp-\yn}h}
\def\dhmun{\frac{\munp-\mun}h}
\def\dhun{\frac{\unp-\un}h}
\def\dhynm{\frac{\yn-\ynm}h}

\def\yk{y^k}
\def\muk{\mu^k}
\def\uk{u^k}

\def\yz{y_0}
\def\mz{m_0}
\def\jz{j_0}
\def\muzl{\mu_0^\lambda}

\def\Lpi{L_\pi}
\def\yu{y^1}
\def\muu{\mu^1}

\Begin{document}

%
\title{Well-posedness and regularity
for a \\generalized fractional Cahn--Hilliard system}
\author{}
\date{}
\maketitle
\Bcenter
\vskip-1cm
{\large\sc Pierluigi Colli$^{(1)}$}\\
{\normalsize e-mail: {\tt pierluigi.colli@unipv.it}}\\[.25cm]
{\large\sc Gianni Gilardi$^{(1)}$}\\
{\normalsize e-mail: {\tt gianni.gilardi@unipv.it}}\\[.25cm]
{\large\sc J\"urgen Sprekels$^{(2)}$}\\
{\normalsize e-mail: {\tt sprekels@wias-berlin.de}}\\[.45cm]
$^{(1)}$
{\small Dipartimento di Matematica ``F. Casorati'', Universit\`a di Pavia}\\
{\small and Research Associate at the IMATI -- C.N.R. Pavia}\\
{\small via Ferrata 5, 27100 Pavia, Italy}\\[.2cm]
$^{(2)}$
{\small Department of Mathematics}\\
{\small Humboldt-Universit\"at zu Berlin}\\
{\small Unter den Linden 6, 10099 Berlin, Germany}\\[2mm]
{\small and}\\[2mm]
{\small Weierstrass Institute for Applied Analysis and Stochastics}\\
{\small Mohrenstrasse 39, 10117 Berlin, Germany}\\[.7cm]
{\it \jupiter{Dedicated to the memory of Professor Emilio Gagliardo in admiration}}
\Ecenter
\Begin{abstract}\noindent
In this paper, we investigate a rather general system of two operator equations
that has the structure of a viscous or nonviscous Cahn--Hilliard 
system in which nonlinearities of
double-well type occur. Standard cases like regular or logarithmic potentials, as well
as non-differentiable potentials involving indicator functions, are admitted. 
The operators appearing in the
system equations are fractional versions of general linear operators \,$A$ and $B$, 
where the latter are densely defined,  unbounded, self-adjoint and monotone in a Hilbert space of functions defined in a smooth domain and have compact resolvents.
In this connection, we remark the fact that our definition of the fractional power of
operators uses the approach via spectral theory. Typical cases are given by standard 
second-order elliptic differential
operators (e.g., the Laplacian) with zero Dirichlet or Neumann boundary 
conditions, but also other cases like fourth-order systems or systems involving
the Stokes operator are covered by the theory. We derive in this paper general
well-posedness and regularity results that extend corresponding results
which are known for either the non-fractional Laplacian with zero Neumann 
boundary condition or the fractional Laplacian with zero \jupiter{Dirichlet condition.  
It turns out that} the first eigenvalue \,$\lambda_1$\, of $\,A\,$ plays an
important und not entirely obvious role: if $\,\lambda_1\,$ is positive, then the
operators $\,A\,$ and $\,B\,$ may be completely unrelated; if, however, $\lambda_1$\, equals
zero, then it must be simple and the corresponding one-dimensional eigenspace has to
consist of the constant functions and to be a subset of the domain of definition
of a certain fractional power of \,$B$. We are able to show general existence,
uniqueness, and regularity results for both these cases, as well as for both the
viscous and the nonviscous system.

\End{abstract}
\salta
\pagestyle{myheadings}
\newcommand\testopari{\sc Colli \ --- \ Gilardi \ --- \ Sprekels}
\newcommand\testodispari{\sc Generalized fractional Cahn--Hilliard system}
\markboth{\testopari}{\testodispari}
\finqui
%

\section{Introduction}
\label{Intro}
\setcounter{equation}{0}

Let $\Omega\subset \erre^3$ denote a bounded, connected and smooth set and
$H$ be a Hilbert space of real-valued functions defined on $\Omega$. We investigate in this paper the abstract evolutionary system
\Bsist
  && \dt y + A^{2r} \mu = 0,
  \label{Iprima}
  \\
  && \tau \dt y + B^{2\sigma} y + f'(y) = \mu + u,
  \label{Iseconda}
  \\
  && y(0) = \yz,
  \label{Icauchy}
\Esist
where $A^{2r}$ and $B^{2\sigma}$, with $r>0$ and $\sigma>0$,
 denote fractional powers of the selfadjoint,
monotone and unbounded linear operators $A$ and $B$, respectively, which are densely defined 
in $H$ and have compact resolvents. The above system can be seen as a generalization
of the famous Cahn--Hilliard system, which models a  
phase separation process taking place in the container $\Omega$ (the list 
\cite{CahH, EZ, NovCoh, FG, Gu, RH, CF2, CGS13} combines basic 
references with some recent contribution on Cahn--Hilliard systems). 
In this case,
one typically has  $A^{2r}=B^{2\sigma}=-\Delta$ with zero Neumann boundary conditions, and the unknown 
functions $\,y\,$ and $\,\mu\,$ stand for the \emph{order parameter} (usually a scaled
density of one of the involved phases) and the \emph{chemical potential}
associated with the phase transition, respectively. Moreover, $\,f\,$ denotes a 
double-well potential.
Typical and physically significant examples for $\,f\,$ 
are the so-called {\em classical regular potential}, the {\em logarithmic double-well potential\/},
and the {\em double obstacle potential\/}, which are given, in this order,~by
\Bsist
  && f_{reg}(r) := \frac 14 \, (r^2-1)^2 \,,
  \quad r \in \erre, 
  \label{regpot}
  \\
  && f_{log}(r) := \bigl( (1+r)\ln (1+r)+(1-r)\ln (1-r) \bigr) - c_1 r^2 \,,
  \quad r \in (-1,1),
  \label{logpot}
  \\[1mm]
  && f_{2obs}(r) := - c_2 r^2 
  \quad \hbox{if $|r|\leq1$}
  \aand
  f_{2obs}(r) := +\infty
  \quad \hbox{if $|r|>1$}.
  \label{obspot}
\Esist
Here, the constants $c_i$ in \eqref{logpot} and \eqref{obspot} satisfy
$c_1>1$ and $c_2>0$, so that $f_{log}$ and $f_{2obs}$ are nonconvex.
In cases like \eqref{obspot}, one has to split $f$ into a nondifferentiable convex part~$\Beta$ 
(the~indicator function of $[-1,1]$, in the present example) and a smooth perturbation~$\Pi$.
Accordingly, one has to replace the derivative of the convex part
by the subdifferential and interpret \eqref{Iseconda} as a differential inclusion
or, equivalently, as a variational inequality
involving \,$\Beta$\, rather than its subdifferential. Actually, we will do
the latter in this paper. We also note that $\tau$ is a nonnegative parameter, 
where for the classical
Cahn--Hilliard system one has $\tau=0$ (the \emph{nonviscous} case); in this paper,
we will handle both the nonviscous case $\tau=0$ and the \emph{viscous} case
$\tau>0$ simultaneously. Of course, better regularity results are to be expected
in the latter case.

Fractional operators are nowadays a very hot topic in \juerg{the} mathematical literature, and it 
occurs that different variants of fractional operators may be considered and 
tackled. Let us perform some review of contributions and results. The paper \cite{Kwa} 
deals with several definitions of the fractional Laplacian (also known as the Riesz 
fractional derivative operator), which is a core 
example of a class of nonlocal pseudodifferential operators appearing in various areas 
of theoretical and applied mathematics. In connection with such fractional operators, 
fractional Sobolev spaces are revisited and discussed in \cite{DPV}. The contributions 
by Servadei and Valdinoci deserve some attention:
in \cite{SV1}, a comparison is made between the spectrum of two different fractional 
Laplacian operators, of which the second one fits in our framework; the paper 
\cite{SV2} discusses the regularity of the weak solution to the fractional Laplace 
equation; the existence of nontrivial solutions for nonlocal semilinear Dirichlet 
problems is established in \cite{SV0}; a fractional counterpart to the well-known 
Brezis--Nirenberg result on the existence of nontrivial solutions to elliptic equations 
with critical nonlinearities is provided in \cite{SV3}.

The paper \cite{AD} presents a construction of harmonic functions on bounded domains 
for the spectral fractional Laplacian operator having a divergent profile at the 
boundary. In the contribution \cite{CT},
a nonlinear pseudodifferential boundary value problem is investigated in a bounded 
domain with homogeneous Dirichlet boundary conditions, where the square root of 
the negative Laplace operator is involved. Regularity results and sharp estimates are 
proved in \cite{CS} for fractional elliptic equations. A nonlocal diffusion operator 
having \juerg{the fractional Laplacian as a special case} is analyzed in \cite{DG} on bounded 
domains, with respect to nonlocal interactions. Fractional Dirichlet and Neumann type 
boundary problems associated with the fractional Laplacian are investigated in  
\cite{Gru1}, by \juerg{demonstrating} regularity properties with a spectral approach; this 
analysis is extended to the fractional heat equation in \cite{Gru2}. 
Obstacle problems for the spectral fractional Laplacian are studied in \cite{MN}. 
By using the Caputo variant of an integral operator with the Riesz kernel, the 
authors of \cite{RS1, RS2} prove regularity up to the boundary for a Dirichlet-type 
boundary value problem and study the extremal solutions by extending some well-known 
results on the extremal solutions when the operator is the Laplacian. \newpier{Some nonlocal problems involving the fractional $p -$Laplacian and nonlinearities at critical growth are examined in \cite{BSY}.}

\juerg{Fractional} porous medium type equations are discussed in \cite{BFV,BSV,BV}. 
The paper \cite{BSV} deals with existence, uniqueness and asymptotic behavior of the 
solutions to an integro-differential equation related to porous medium equations \juerg{in}  
bounded domains; the problem does not have a separate boundary condition, since zero 
boundary data are implicitly assumed in the definition of the operator. A priori 
estimates for positive solutions of a porous medium equation are shown in \cite{BV}, 
where the spectral fractional Laplacian with zero Dirichlet boundary data is 
considered; it turns out that the results are influenced by the first eigenvalue and 
eigenfunction. A quantitative study of nonnegative solutions of the same equation  is provided in 
\cite{BFV}, where the regularity theory is addressed: decay and positivity, Harnack 
inequalities, interior and boundary regularity, and asymptotic behavior are investigated. \newpier{Also fractional Schr\"odinger equations are receiving a good deal of attention, see, e.g., \cite{BSX} and references therein.} 

\juerg{There exist} already quite a number of contributions 
dealing with nonlocal variants of the Cahn--Hilliard system. In \cite{ABG}, the problem 
of well-posedness for a nonlocal Cahn--Hilliard equation is established by 
interpreting the problem as a Lipschitz perturbation of a maximal monotone operator in 
a suitable Hilbert space. A fractional variant of the Cahn--Hilliard equation settled 
in a bounded domain and complemented with homogeneous Dirichlet boundary conditions of 
solid type is introduced in \cite{AkSS1}: existence and uniqueness of weak solutions 
to the related initial-boundary value problem are proved and some significant singular 
limits are investigated as the order of either of the fractional Laplacians appearing 
in the system \juerg{approaches zero}. Moreover, in the recent paper \cite{AkSS2}, for fixed 
orders of the operators the convergence as time goes to infinity of each solution to a 
(single) equilibrium is proved. In \cite{AM}, the authors derive a fractional Cahn--Hilliard 
equation by considering a gradient flow in the negative order Sobolev space $H^{-\alpha}$, 
$\alpha\in [0,1]$, 
where the choice $\alpha=1$ corresponds to the classical Cahn--Hilliard equation, while 
the choice $\alpha=0$ recovers the Allen--Cahn equation; existence and stability estimates 
are derived in the case where the nonlinearity is a quartic polynomial, as in \eqref{regpot}. 
The paper \cite{GGG} addresses the nonlocal Cahn--Hilliard equation with a singular potential 
and a constant mobility: among a class of results, in particular the authors can establish 
the validity of the strict separation property in \juerg{two dimensions}. Another interesting analysis 
of a nonstandard and nonlocal Cahn--Hilliard system can be found in \cite{CGS6}. Next, in\cite{GalEJAM}
a nonlocal version of the Cahn--Hilliard equation characterized by the presence of a fractional 
diffusion operator, and which is subject to fractional dynamic boundary conditions, is studied. 
The articles \cite{GalDCDS,GalAIHP} 
\juerg{treat} a doubly nonlocal Cahn--Hilliard equation with special kernels in the operators: 
well-posedness results, along with regularity, long-time behavior, and 
global attractors, are investigated in connection with the interaction between the two levels 
of nonlocality in the operators. \jupiter{In particular, the paper which is closest to our results and partially covers some of them, is \cite{GalDCDS}. There, also results for the logarithmic potential are provided.}

In our approach,  \juerg{which we develop} in the subsequent sections, we work with fractional 
operators defined via spectral theory. This position \juerg{enables} us to deal with 
powers of a second-order elliptic operator with either Dirichlet or Neumann or 
Robin boundary conditions, allowing us a wide setting in this respect. Moreover, 
other operators, \juerg{such as} fourth-order ones or systems involving
the Stokes operator, can \juerg{be covered by} the theory. 

The aim of the present paper is to prove
general well-posedness and regularity theorems that extend the corresponding 
results known for either the non-fractional Laplacian with zero Neumann 
boundary condition or the fractional Laplacian with zero Dirichlet condition 
(cf.~\cite{AM,AkSS1}). In the development of the theory, one realizes that 
 the first eigenvalue \,$\lambda_1$\, of $\,A\,$ plays an
important und not entirely obvious role. Indeed, it turns out that 
 if $\,\lambda_1\,$ is positive, then the
operators $\,A\,$ and $\,B\,$ may be completely unrelated. 
On the other hand, in the case when $\lambda_1 = 0$, \juerg{then} 
we have to assume that $\lambda_1$ is \juerg{a simple eigenvalue and that}
the corresponding one-dimensional eigenspace
consists of constant functions, on which 
the proper fractional power of \,$B$ should operate. This set of assumptions looks \juerg{like} 
a heavy restriction, but let us notice that the framework is strongly related 
to the structure of the Cahn--Hilliard system with the natural Neumann 
homogeneous boundary conditions (that exactly imply conservation of mass). In conclusion, 
it \juerg{will turn} out that we are able to show well-posedness and regularity results for both the abovementioned situations, as well as for both the viscous and the nonviscous system, under very general  
assumptions for the convex parts  of the potential $f$ (see~\eqref{regpot}--\eqref{obspot}). 

Here is a brief outline of the paper. Section~\ref{STATEMENT} contains a precise 
statement of the problem along with assumptions and main results; some remarks 
commenting the results and introducing examples of operators are also included. 
Section~\ref{AUX} is intended to present some auxiliary material about relations 
among the involved spaces and properties of the operators; all this turns to be a useful 
toolbox for the \juerg{following analysis}. Section~\ref{UNIQUENESS} deals with the continuous 
dependence of the solution on the data, while Section~\ref{APPROX} introduces 
an approximating problem
based on the Moreau--Yosida regularizations of the convex functions and on an 
implicit time discretization of the system, which is fully discussed concerning 
existence of the discrete solution and uniform a priori estimates \juerg{for} it. 
Section~\ref{EXISTENCE} brings the existence proof, which is carried out 
by taking the limits with respect to the approximation parameters. Finally, 
Section~\ref{REGULARITY} is devoted to show the proper estimates ensuring 
the regularity properties for the solution.


\section{Statement of the problem and results}
\label{STATEMENT}
\setcounter{equation}{0}

In this section, we state precise assumptions and notations and present our results.
First of all, the set $\Omega\subset\erre^3$ is assumed to be bounded, connected and 
smooth,
with outward unit normal vector field $\,\nu\,$ on $\Gamma:=\partial\Omega$.
Moreover, $\dn$ stands for the corresponding normal derivative.
We use the notation
\Beq
  H := \Ldue
  \label{defH}
\Eeq
and denote by $\norma\cpto$ and $(\cpto,\cpto)$ the standard norm and inner product of~$H$.
Now, we start introducing our assumptions.
We first postulate that
\Bsist
  && A:D(A)\subset H\to H
  \aand
  B:D(B)\subset H\to H
  \quad \hbox{are}
  \non
  \\
  && \hbox{unbounded monotone selfadjoint linear operators with compact resolvents.} 
  \qquad
  \label{hpAB} 
\Esist
This assumption implies that there are sequences 
$\{\lambda_j\}$ and $\{\lambda'_j\}$ of eigenvalues
and orthonormal sequences $\{e_j\}$ and $\{e'_j\}$ of corresponding eigenvectors,
that~is,
\Beq
  A e_j = \lambda_j e_j, \quad
  B e'_j = \lambda'_j e'_j
  \aand
  (e_i,e_j) = (e'_i,e'_j) = \delta_{ij}
  \quad \hbox{for $i,j=1,2,\dots$}
  \label{eigen}
\Eeq
\pier{with $\delta_{ij}$ denoting the Kronecker index,} such that
\begin{align}
  & 0 \leq \lambda_1 \leq \lambda_2 \leq \dots
  \aand
  0 \leq \lambda'_1 \leq \lambda'_2 \leq \dots
  \quad \hbox{with} \quad
  \lim_{j\to\infty} \lambda_j
  = \lim_{j\to\infty} \lambda'_j
  = + \infty,
  \label{eigenvalues}
  \\[1mm]
  & \hbox{$\{e_j\}$ and $\{e'_j\}$ are complete systems in $H$}.
  \label{complete}
\end{align}
The above assumptions on $A$ and $B$ allow us to define the powers of $A$ and $B$ 
for an arbitrary positive real exponent.
As far as the first operator is concerned, we have
\Bsist
  && \VA r := D(A^r)
  = \Bigl\{ v\in H:\ \somma j1\infty |\lambda_j^r (v,e_j)|^2 < +\infty \Bigr\}
  \aand
  \label{defdomAr}
  \\
  && A^r v = \somma j1\infty \lambda_j^r (v,e_j) e_j
  \quad \hbox{for $v\in\VA r$},
  \label{defAr}
\Esist
the series being convergent in the strong topology of~$H$,
due to the properties \eqref{defdomAr} of the coefficients.
In principle, we endow $\VA r$ with the (graph) norm and inner product
\Beq
  \norma v_{gr,A,r}^2 := (v,v)_{gr,A,r}
  \aand
  (v,w)_{gr,A,r} := (v,w) + (A^r v , A^r w)
  \quad \hbox{for $v,w\in\VA r$}.
  \label{defnormagrAr}
\Eeq
This makes $\VA r$ a Hilbert space.
However, we can choose any equivalent Hilbert norm.
Later on, we actually will do that.
In the same way, starting from \accorpa{hpAB}{complete} for~$B$,
we can define the power $B^\sigma$ of $B$ for every $\sigma>0$.
We therefore set
\Bsist
  && \VB\sigma := D(B^\sigma),
  \quad \hbox{with the norm $\norma\cpto_{B,\sigma}$ associated to the inner product}
  \label{defBs}
  \non
  \\
  && (v,w)_{B,\sigma} := (v,w) + (B^\sigma v,B^\sigma w)
  \quad \hbox{for $v,w\in \VB\sigma$}.
  \label{defprodBs}
\Esist
If $r_i$ and $\sigma_i$ are arbitrary positive exponents,
it is clear that
\Bsist
  && (A^{r_1+r_2} v,w)
  = (A^{r_1} v, A^{r_2} w)
  \quad \hbox{for every $v\in\VA{r_1+r_2}$ and $w\in\VA{r_2}$},
  \label{propA}
  \\[1mm]
  && (B^{\sigma_1+\sigma_2} v,w)
  = (B^{\sigma_1} v, B^{\sigma_2} w)
  \quad \hbox{for every  $v\in\VB{\sigma_1+\sigma_2}$ and $w\in\VB{\sigma_2}$}.
  \label{propB}
\Esist
\noindent
From now on, we assume:
\Beq
  \hbox{$r$ and $\sigma$ are fixed positive real numbers.}
  \label{hprs}
\Eeq
\gianni{Accordingly, we introduce a space with a negative exponent}.
We set
\Beq
  \VA{-r} := (\VA r)^* 
  \quad \hbox{for $r>0$},
  \label{defVAneg}
\Eeq
and use the symbol $\<\cpto,\cpto>_{A,r}$ for the duality pairing
between $\VA{-r}$ and~$\VA r$.
Moreover, we identify $H$ with a subspace of $\VA{-r}$
in the usual way, i.e., such that
\Beq
  \< v,w >_{A,r} = (v,w)
  \quad \hbox{for every $v\in H$ and $w\in\VA r$}.
  \label{identification}
\Eeq
Next, we make the following assumption:
\Bsist
  && \hbox{Either} \quad
  \lambda_1 > 0 
  \quad \hbox{or} \quad
  \hbox{$0=\lambda_1<\lambda_2$ and $e_1$ is a constant.}
  \label{hpsimple}
  \\
  && \hbox{\gianni{If\quad $\lambda_1=0$,\quad}
  the constant functions belong to $\VB\sigma$}.
  \label{hpVB}
\Esist

\Brem
\label{Remhpsimple}
Let us comment on the assumptions~\eqref{hpsimple}.
The meaning of the first case is clear, and such a condition is satisfied
by the more usual elliptic operators with zero Dirichlet boundary conditions
(however, also mixed boundary conditions could be considered, 
with proper definitions of the domains of the operators),
for instance:
$i)$~$A$ is the Laplace operator $-\Delta$ with domain $D(-\Delta)=\Hdue\cap\Hunoz$;
$ii)$~$A$ is the bi-harmonic operator $\Delta^2$ with domain: 
$D(\Delta^2)=\Hx4\cap H^2_0(\Omega)$.
The second case of \eqref{hpsimple},
where the strict inequality means that the first eigenvalue $\lambda_1=0$ is simple,
happens in the following important situations:
$i)$~$A$ is the Laplace operator $-\Delta$ with zero Neumann boundary conditions,
which corresponds to the choice $D(-\Delta)=\{v\in\Hdue:\ \dn v=0\}$;
$ii)$~$A$ is the bi-harmonic operator $\Delta^2$ with the boundary conditions
corresponding to the following choice of the domain: 
$D(\Delta^2)=\{v\in\Hx4:\ \dn v=\dn\Delta v=0\}$.
Indeed, $\Omega$~is assumed to be bounded, smooth and connected.
\Erem

\Brem
\label{RemhpVB}
\gianni{We point out that \eqref{hpVB} is the only condition that involves both operators $A$ and~$B$,
i.e., if $\lambda_1>0$, these operators are completely unrelated.
However, we notice that the assumption on the constant functions is rather mild.
Indeed, it} 
holds for many operators
whose domain involves Neumann boundary conditions.
This is the case, for instance,
if $B$ is the Laplace operator with domain $D(-\Delta)=\{v\in\Hdue:\ \dn v=0\}$.
On the contrary, if $B=-\Delta$ with domain $D(-\Delta):=\Hdue\cap\Hunoz$,
then $D(B)$ does not contain any nonzero constant functions.
However, $\VB\sigma$ does contain  every constant function
provided that $\sigma\in(0,1/4)$,
since it coincides with the usual Sobolev-Slobodeckij space~$\Hx{2\sigma}$.
Indeed, the spaces $\VA r$ and $\VB\sigma$ can be seen in the framework of interpolation theory.
However, we prefer to avoid this check and deduce all the results we need 
from our definitions.
\Erem

\Brem
\label{Stokes}
We have chosen to take $H:=\Ldue$ once and for all, for simplicity.
However, it is clear that our assumptions are rather close to an abstract situation
and can be adapted to other choices of the space~$H$ as well.
For instance, one could deal with the Stokes operator with Dirichlet boundary conditions,
by taking for $H$ the space of vector-valued functions $v\in(\Ldue)^3$ satisfying
$\div v=0$ in the sense of distributions and defining the operator $A$ as follows:
\pier{an element $v\in H$ belongs to $D(A)$ if and only if $v \in (\Hunoz)^3$ and 
$\Delta v:=(\Delta v_i) \in (\Ldue)^3$; for $v\in D(A)$, $Av$ is the $L^2$-projection on $H$ of $- \Delta v$.}
In this case, the first assumption of \eqref{hpsimple} is satisfied.
Of course, the hypotheses on the structure of the nonlinear terms to be 
introduced below would have to be adapted to this new situation.
\Erem

We use assumption \eqref{hpsimple} to define a different Hilbert norm on~$\VA r$.
We set, for $v\in\VA r$,
\Beq
  \norma v_{A,r}^2 := \left\{ 
  \begin{aligned}
  & \norma{A^r v}^2
  = \somma j1\infty |\lambda_j^r (v,e_j)|^2
  \qquad \hbox{if $\lambda_1>0$,}
  \\
  & |(v,e_1)|^2 + \norma{A^r v}^2
  = |(v,e_1)|^2 + \somma j2\infty |\lambda_j^r (v,e_j)|^2
  \qquad \hbox{if $\lambda_1=0$}.
  \end{aligned}
  \right.
  \label{defnormaAr}
\Eeq
In the next section we will show that this norm is equivalent 
to the graph norm defined in~\eqref{defnormagrAr},
and we always will use the norm \eqref{defnormaAr} rather than~\eqref{defnormagrAr}.
Of course, we will also use the corresponding inner product in $\VA r$ and norm in~$\VA{-r}$.
They are given~by
\Bsist
  && (v,w)_{A,r}
  = (A^r v,A^r w)
  \quad \hbox{or} \quad
  (v,w)_{A,r}
  = (v,e_1)(w,e_1) + (A^r v,A^r w),
  \qquad
  \non
  \\
  && \quad \hbox{depending on whether $\lambda_1>0$ or $\lambda_1=0$,\quad 
	for $v,w\in\VA r$,}
  \label{defprodAr}
  \\[1mm]
  && \hbox{$\norma\cpto_{A,-r}$ is the dual norm of $\norma\cpto_{A,r}$}.
  \label{defnormaVAneg}
\Esist

\Brem
\label{RemnormaVAr}
We notice that in the case $\lambda_1=0$ of~\eqref{hpsimple}
the constant value of $e_1$
is equal to one of the numbers $\pm|\Omega|^{-1/2}$,
where $|\Omega|$ is the volume of~$\Omega$.
It follows for every $v\in H$ that the first term $(v,e_1)e_1$ of the Fourier series of $v$
is the constant function whose value~is the mean value of~$v$,~i.e.,
\Beq
  \mean v := \frac 1 {|\Omega|} \iO v\,,
  \label{defmean}
\Eeq
and that the first terms of the sums appearing in \eqref{defnormaAr} and \eqref{defprodAr}
are given~by
\Bsist
  && |(v,e_1)|^2 = |\Omega| \, (\mean v)^2
  \quad \hbox{for every $v\in H$},
  \non
  \\[1mm]
  && (v,e_1) (w,e_1) = |\Omega| \, (\mean v) (\mean w)
  \quad \hbox{for every $v,w\in H$}.
  \non
\Esist
\Erem

For the other ingredients of our system, we postulate the following properties:
\begin{align}
  & \hbox{$\tau$ is a nonnegative real number.}
  \label{hptau}
  \\
  & \Beta : \erre \to [0,+\infty]
  \quad \hbox{is convex, proper and l.s.c.\ with} \quad
  \Beta(0) = 0.
  \label{hpBeta}
  \\
  \separa
  & \Pi : \erre \to \erre
  \quad \hbox{is of class $C^1$ with a \Lip\ continuous first derivative.}
  \label{hpPi}
  \\
  & \mbox{{It holds }}\,\,\displaystyle \liminf_{|s|\nearrow+\infty} \displaystyle \frac {\Beta(s)+\Pi(s)} {s^2} > 0.
  \label{hpcoerc}
\end{align}
We can suppose that $\tau\leq1$ without loss of generality.
We remark that the assumptions \accorpa{hpBeta}{hpcoerc} are 
fulfilled by all of the important potentials \accorpa{regpot}{obspot}.
We set, for convenience,
\Beq
  \beta := \partial\Beta , \quad
  \pi := \Pi' , \quad
  \Lpi = \hbox{the \Lip\ constant of $\pi$,}
  \aand
  \Lpi' := \Lpi + 1 \,.
  \label{defbetapi}  
\Eeq
Moreover, we term $D(\Beta)$ and $D(\beta)$ the effective domains of $\Beta$ and~$\beta$, respectively,
and, for $r\in D(\beta)$,
we use the symbol $\betaz(r)$ for the element of $\beta(r)$ having minimum modulus.
Notice that $\beta$ is a maximal monotone graph in $\erre\times\erre$.

At this point, we can state the problem {under investigation}.
On account of \accorpa{propA}{propB}, we give a weak formulation of the equations \accorpa{Iprima}{Iseconda}.
Moreover, we present \eqref{Iseconda} as a variational inequality.
For the data, we make the following assumptions:
\Bsist
  && u \in \H1H 
  \label{hpu}
  \\
  && \yz \in \VB\sigma
  \aand
  \Beta(\yz) \in \Luno 
  \label{hpyz}
  \\
  && \gianni{\hbox{if $\lambda_1=0$,}} \quad
  \mz := \mean\yz
  \quad \hbox{belongs to the interior of $D(\beta)$} .
  \label{hpmz}
\Esist
\Accorpa\HPdati hpu hpmz
\gianni{Notice that no condition on $\mz$ is required if $\lambda_1>0$.}
Then, we set
\Beq
  Q := \Omega \times (0,T)
  \label{defQ}
\Eeq
and look for a pair $(y,\mu)$ satisfying
\Bsist
  && y \in \H1{\VA{-r}} \cap \L\infty{\VB\sigma} 
  \aand
  \tau \dt y \in \L2H,
  \qquad
  \label{regy}
  \\
  && \mu \in \L2{\VA r},
  \label{regmu}
  \\
  && \Beta(y) \in \LQ1,
  \label{regBetay}
\Esist
\Accorpa\Regsoluz regy regBetay
and solving the system
\Bsist
  && \< \dt y(t) , v >_{A,r}
  + ( A^r \mu(t) , A^r v )
  = 0
  \quad \hbox{for every $v\in\VA r$\, and \,a.e. \,}t\in (0,T),
  \qquad
  \label{prima}
  \\[2mm]
  && \previousgianni{\bigl( \tau \dt y(t) , y(t) - v \bigr)}
  + \bigl( B^\sigma y(t) , B^\sigma( y(t)-v) \bigr)
  \non
  \\
  && \quad {}
  + \iO \Beta(y(t))
  + \bigl( \pi(y(t)) - u(t) ,  y(t)-v \bigr)
  \leq \bigl( \mu(t) ,  y(t)-v \bigr)
  + \iO \Beta(v)
  \non
  \\
  && \quad \hbox{for every $v\in\VB\sigma$\, and \, a.e. \,} t\in (0,T),
  \label{seconda}
  \\[0.2cm]
  && y(0) = \yz \,.
  \label{cauchy}
\Esist
\Accorpa\Pbl prima cauchy
\pier{Of course, it is understood that 
$$ \iO \Beta(v) = +\infty  \quad \hbox{whenever} \quad \Beta (v) \notin \Luno.$$
A similar agreement also holds for \juerg{integrals of the type} $\intQ \Beta(v)$ whenever 
$v\in L^2(Q)$ but $\Beta(v) \notin L^1(Q)$.}

\pier{Now, let us notice} that \eqref{seconda} is equivalent to its time-integrated variant, that is,
\begin{align}
  & \previousgianni{\ioT \bigl(\tau \dt y(t) , y(t) - v \bigr) \, dt }
  + {\ioT \bigl( B^\sigma y(t) , B^\sigma( y(t)-v(t)) \bigr)\,dt }
  \non
  \\ 
  & \quad {{}
  + \intQ \Beta(y) 
  + \ioT \bigl( \pi(y(t)) - u(t) ,  y(t)-v(t) \bigr)\,dt} 
  \non
  \\
  & \leq {\ioT \bigl( \mu(t) ,  y(t)-v(t) \bigr)\,dt
  + \intQ \Beta(v)
  \quad \hbox{for every $v\in\L2{\VB\sigma}$}}.
  \qquad
  \label{intseconda}
\end{align}
We also remark that, if $\lambda_1=0$, then $A^r(1)=0$ by \eqref{hpsimple},
so that \eqref{prima} implies that
\Beq
  \frac d{dt} \iO y(t) = 0
  \quad \aat, \quad
  \hbox{i.e.,} \quad
  \mean y(t) = \mz
  \quad \hbox{for every $t\in[0,T]$}.
  \label{conservation}
\Eeq
Finally, \pier{let us note that if $\lambda_1=0$\juerg{, then} the condition \eqref{hpmz} on $\mz$} 
ensures the existence of some $\delta_0>0$ satisfying
\Beq
  [\mz-\delta_0,\mz+\delta_0] \subset D(\beta).
  \label{defdeltaz}
\Eeq
\pier{\Brem
\label{Signif}
According to the definition of subdifferential (cf., e.g., \cite{Brezis} or \cite{Barbu}), the precise meaning of the inequality \eqref{seconda} is that there exists some element $\chi \in \L2{(V_B^\sigma)^*}$ such that 
$$ \chi:= \mu - \tau \dt y -  B^{2\sigma} y -  \pi(y) + u \in \partial \Phi  (y)   \quad \hbox{a.e. in } (0,T) ,$$
where $\partial \Phi $ is the subdifferential of the convex function 
$\Phi : V_B^\sigma \to [0, +\infty]$ defined by 
$$
\Phi (v) :=  \iO \Beta(v) \quad \hbox{if} \ \, \Beta(v) \in \Luno, \quad \Phi (v) := +\infty \quad \hbox{otherwise,}
$$
and actually the subdifferential $\partial \Phi $ is a maximal monotone operator 
from $V_B^\sigma$ to $(V_B^\sigma)^*$. In this sense, \eqref{seconda} turns out to be 
a slight generalization of \eqref{Iseconda}. 
\Erem}

Here is our well-posedness and continuous dependence result.

\Bthm
\label{Wellposedness}
Let the assumptions \eqref{hpAB}, \gianni{\eqref{hprs}}, \accorpa{hpsimple}{hpVB} and 
\eqref{hptau}--{\eqref{hpcoerc}} 
on the structure of the system,
and \HPdati\ on the data, be fulfilled.
Then there exists \jupiter{at least one} pair $(y,\mu)$ satisfying \Regsoluz\ and solving  problem 
\Pbl. \jupiter{Besides, this solution $(y,\mu)$} satisfies the estimate
\Beq
  \norma y_{\H1{\VA{-r}} \cap \L\infty{\VB\sigma}}
  + \norma\mu_{\L2{\VA r}}
  + \norma{\Beta(y)}_{\LQ1}
  + \norma{\tau^{1/2}\dt y}_{\L2H}
  \leq K_1,
  \label{stimasoluz}
\Eeq
with a constant $K_1$ that depends only 
on the structure of the system, the norms of the data corresponding to \accorpa{hpu}{hpyz},
the width $\delta_0$ satisfying~\eqref{defdeltaz} \gianni{if $\lambda_1=0$}, and~$T$.
Moreover, if $u_i$, $i=1,2$, are two choices of~$u$ and $(y_i,\mu_i)$ \jupiter{are corresponding} solutions, then
we~have
\Beq
  \norma{y_1-y_2}_{\pier{\L\infty{\VA{-r}}{}}\cap\L2{\VB\sigma}}
  + \norma{\tau^{1/2}(y_1-y_2)}_{\L\infty H}
  \leq K_2 \norma{u_1-u_2}_{\L2H},
  \label{contdep}
\Eeq
with a constant $K_2$ that depends only on the operators $A^r$ and $B^\sigma$,
the \Lip\ constant~$\Lpi$, and~$T$. \jupiter{In particular, the first component $y$ of any solution $(y,\mu)$ is uniquely determined. In addition, $A^r(\mu)$ is uniquely determined as well. Finally, also the component $\mu$ is uniquely determined if $\lambda_1 >0$.}
\Ethm

\Brem
\label{Remcontdep}
More generally, we could take two different initial values $y_{0,1}$ and $y_{0,2}$,
by assuming that they have the same mean value if $\lambda_1=0$.
Then, the \rhs\ of \eqref{contdep} has to be modified by adding two contributions 
involving $d_0:=y_{0,1}-y_{0,2}$,
which are proportional to 
$\,\norma{d_0}_{A,-r}$\, and to \,$\tau^{1/2}\norma{d_0}$.
\Erem

Under additional assumptions on the data, we have stronger regularity results
in both the viscous and nonviscous cases.
Namely, we also assume that
either $\tau>0$ and
\Beq
  \yz \in \VB{2\sigma}
  \aand
  \betaz(\yz) \in H
  \label{hpyzreg}
\Eeq
or $\tau=0$ and
\Bsist
  && \yz \in \VB{2\sigma} 
  \aand
  \norma{\muzl(t)}_{A,r} \leq M_0,
  \quad \hbox{where} \quad
  \label{hpdatireg}
  \\
  && \muzl(t) := B^{2\sigma} \yz + (\betal+\pi)(\yz) - u(t), 
  \label{defmul}
\Esist
for some constant $M_0$ and every sufficiently small $\lambda>0$ and $t>0$, 
$\betal$ being the Yosida approximation of $\beta$ at the level~$\lambda$
(see, e.g., \cite[p.~28]{Brezis}).
More precisely, it is assumed that the element~$\muzl(t)$
(which is well defined by \eqref{defmul} due to the first assumption on~$\yz$) 
belongs to $\VA r$ and satisfies the above estimate.
Of course, this assumption is very restrictive.
However, we can give sufficient conditions for~it.
{One possibility is to assume that each} of the four contributions to the \rhs\ 
of \eqref{defmul}
satisfies bounds like~\eqref{hpdatireg}, separately, 
and that $A^r$ is a local operator in order to deal with the term~$\betal(\yz)$.
For instance, if $A^r$ is the Laplace operator with Dirichlet boundary conditions
and $\beta$ is single-valued and smooth in the interior of its domain,
then one can assume that
$\yz\in\Hdue\cap\Hunoz$ and that $\min\yz>\inf D(\beta)$ and $\max\yz<\sup D(\beta)$.
These assumptions keep $\betal(\yz)$ bounded in~$\Hdue$, indeed.

\Bthm
\label{Regularity}
In addition to the assumptions of Theorem~\ref{Wellposedness},
suppose that either $\tau>0$ and \eqref{hpyzreg}
or $\tau=0$ and \accorpa{hpdatireg}{defmul} are fulfilled.
Then \jupiter{the solution $(y,\mu)$ established in Theorem~\ref{Wellposedness}}
also satisfies the regularity properties
\Bsist
  && \dt y \in \L\infty{\VA{-r}} \cap \L2{\VB\sigma} 
  \aand
  \mu \in \L\infty{\VA r}
  \quad \hbox{if $\tau\geq0$},
  \qquad
  \label{regsoluzbis}
  \\
  && \dt y \in \L\infty H
  \aand
  \mu \in \L\infty{\VA{2r}}
  \quad \hbox{if $\tau>0$},
  \label{regsoluzter}
\Esist
as well as the estimate
\Bsist
  && \norma{\dt y}_{\L\infty{\VA{-r}}\cap\L2{\VB\sigma}}
  + \norma\mu_{\L\infty{\VA r}}
  \non
  \\
  && \quad {}
  + \norma{\tau^{1/2}\dt y}_{\L\infty H}
  + \norma{\tau^{1/2}\mu}_{\L\infty{\VA{2r}}}
  \leq K_3,
  \label{stimareg}
\Esist
with a constant $K_3$ that depends only 
on the structure of the system, the norms of the data,
the width $\delta_0$ satisfying~\eqref{defdeltaz} \gianni{if $\lambda_1=0$}, 
the constant $M_0$ satisfying~\eqref{hpdatireg} if $\tau=0$, and~$T$.
\Ethm

The remainder of the paper is organized as follows. 
The next section collects some notations and tools that will prove \pier{to be} useful in the sequel.
\jupiter{The continuous dependence result} is proved in Section~\ref{UNIQUENESS},
while the existence of a solution and its regularity are proved
in the last two Sections~\ref{EXISTENCE} and~\ref{REGULARITY} and are prepared by the study 
of the approximating problem introduced in Section~\ref{APPROX}.


\section{Auxiliary material}
\label{AUX}
\setcounter{equation}{0}

\pier{Here,} we add some comments on the spaces defined in the previous \pier{section}.
Moreover, we introduce some new spaces and operators, 
as well as some notations and properties concerning interpolating functions.
First of all, we stress the following facts:
\Bsist
  && \hbox{The embeddings $\VA{r_2} \subset \VA{r_1} \subset H$ are dense and compact for $0<r_1<r_2$}.
  \label{compembA}
  \\
  && \hbox{The embeddings $H \subset \VA{-r_1} \subset \VA{-r_2}$ are dense and compact for $0<r_1<r_2$}.
  \qquad
  \label{compembAneg}
  \\
  && \hbox{The embeddings $\VB{\sigma_2} \subset \VB{\sigma_1} \subset H$ are dense and compact for $0<\sigma_1<\sigma_2$}.
  \label{compembB}
\Esist
Let us comment on just the first embedding of \eqref{compembA},
since the second one and \eqref{compembB} are similar 
and \eqref{compembAneg} follows as a consequence of~\eqref{compembA}.
The density is clear.
For compactness, notice that
$\lim_{j\to\infty}\lambda_j^{r_1-r_2}=0$,
so that the mapping  that to each \,$\{c_j\}\in\ell^2$\,
associates \,$\{\lambda_j^{r_1-r_2}c_j\}$\,
is compact from $\ell^2$ into itself.

From the continuous embedding $H\subset\VA{-r}$ 
and the compact embedding $\VB\sigma\subset H$ given by~\accorpa{compembAneg}{compembB},
it follows that, for every $\delta>0$, there exists a constant $c_\delta$ such that
\Beq
  \norma v^2
  \leq \delta \, \norma{B^\sigma v}^2 
  + c_\delta \norma v_{\VA{-r}}^2
  \quad \hbox{for every $v\in\VB\sigma$}.
  \label{compineq}
\Eeq

\Bprop
\label{Equivnorms}
The norms \eqref{defnormagrAr} and~\eqref{defnormaAr} on $\VA r$ are equivalent.
\Eprop

\Bdim
Take any $v\in H$.
Then, $v$~can be represented in the form
\Beq
  v = \somma j1\infty c_j e_j
  \quad \hbox{in $H$,}
  \quad \hbox{where $c_j:=(v,e_j)$ \,\,for all $j\in\enne$, }
  \non
\Eeq
and where the sequence $\{c_j\}_{j\geq1}$ belongs to~$\ell^2$.
On the other hand, by the definition of~$\VA r$, $A^r$ and~$\norma\cpto_{gr,A,r}$, 
we have, for every $v\in H$,
\Beq
  v\in\VA r
  \quad \hbox{if and only if} \quad
  \somma j1\infty |\lambda_j^r c_j|^2 < + \infty ,
  \aand
  \norma v_{gr,A,r}^2
  = \norma v^2 + \somma j1\infty |\lambda_j^r c_j|^2 .
  \non
\Eeq
Therefore, by recalling~\eqref{defnormaAr}, we conclude that
\Beq
  \norma v_{A,r} \leq \norma v_{gr,A,r}.
  \non
\Eeq
Now, suppose that $\lambda_1>0$.
Then we have
\Beq
  \lambda_1^{2r} \norma v^2 
  = \lambda_1^{2r} \somma j1\infty |c_j|^2
  \leq \somma j1\infty |\lambda_j^r c_j|^2
  = \norma{A^r v}^2,
  \non
\Eeq
since $\lambda_j\geq\lambda_1$ for every~$j$, 
whence immediately 
\Beq
  \norma v_{gr,A,r}^2
  = \norma v^2 + \norma{A^r v}^2
  \leq \Bigl( \frac 1 {\lambda_1^{2r}} + 1 \Bigr) \norma{A^r v}^2
  = \Bigl( \frac 1 {\lambda_1^{2r}} + 1 \Bigr) \norma v_{A,r}^2 \,.
  \non
\Eeq
If, instead, $\lambda_1=0$, then we recall that $A^r(1)=0$ and thus
\Bsist
  && \norma v_{gr,A,r}^2
  \leq 2 \norma{\mean v}_{gr,A,r}^2
  + 2 \norma{v-\mean v}_{gr,A,r}^2
  \non
  \\[1mm]
  && = 2 \norma{\mean v}^2
  + 2 \bigl( \norma{v-\mean v}^2 + \norma{A^r(v-\mean v)}^2 \bigr)
  \non
  \\[1mm]
  && = 2 |\Omega| \, |\mean v|^2
  + 2  \norma{v-\mean v}^2
  + 2 \norma{A^r v}^2 \,.
  \non
\Esist
Thus, on accout of Remark~\ref{RemnormaVAr}, the desired inequality follows if we prove that,
for some constant $\hat c>0$,  it holds the Poincar\'e type inequality
\Beq
  \norma v \leq \hat c \, \norma{A^r v}
  \quad \hbox{for every $v\in\VA r$ with $\mean v=0$}.
  \label{poincare} 
\Eeq
This is an easy consequence of the compact embedding
$\VA r\subset H$ (see~\eqref{compembA}).
However, we prove it for the reader's convenience.
By contradiction, there exists a sequence $\{v_n\}$ in $\VA r$ satisfying
\Beq
  \norma{v_n} > n \norma{A^r v_n}
  \aand 
  \mean v_n = 0
  \quad \hbox{for every $n\geq1$}.
  \non
\Eeq
Clearly, we have that $v_n\not=0$, so that
we can define $w_n:=v_n/\norma{v_n}$.
Then, $\norma{w_n}=1$, $\norma{A^r w_n}<1/n$ and $\mean w_n=0$ for every~$n$.
In particular, $\{w_n\}$ is bounded in $\VA r$, whence we have
\Beq
  w_{n_k} \to w
  \quad \hbox{weakly in $\VA r$}
  \non
\Eeq
for some subsequence and some $w\in\VA r$.
By the compact embedding $\VA r\subset H$,
we infer that $w_{n_k}$ converges to $w$ strongly in~$H$,
whence $\norma w=1$ and $\mean w=0$.
On the other hand, we also have that $A^r w=0$ since $\norma{A^r w_n}<1/n$ for every~$n$.
Therefore, $w$~is a constant.
Hence, the above conclusions $\norma w=1$ and $\mean w=0$ yield a contradiction.
\Edim

At this point, we introduce the Riesz isomorphism $\calR_r:\VA r\to\VA{-r}$
associated with the inner product~\eqref{defprodAr}, which acts as follows:
\Beq
  \< \calR_r v , w >_{A,r}
  = (v,w)_{A,r}
  \quad \hbox{for every $v,w\in\VA r$}.
  \label{riesz}
\Eeq
Moreover, we~set
\begin{align}
  & \Vz r := \VA r
  \aand
  \Vz{-r} := \VA{-r}
  \quad \hbox{if $\lambda_1>0$},
  \non
  \\[1mm]
  & \Vz r := \{v\in \VA r :\ \mean v=0\}
  \aand
  \Vz{-r} := \{v \in \VA{-r} :\ \<v,1>_{A,r}=0 \}
  \quad \hbox{if $\lambda_1=0$} \,.
  \label{defVrpos}
\end{align}

\Bprop
\label{ExtensionR}
The Riesz isomorphism $\calR_r$ maps $\Vz r$ onto $\Vz{-r}$.
Moreover, $\calR_r$ extends to $\Vz r$ the restriction of $A^{2r}$ to $\Vz{2r}$.
\Eprop

\Bdim
Let us deal with the first assertion.
If $\lambda_1>0$, there is nothing to prove.
Thus, assume that $\lambda_1=0$. 
Then, on account of Remark~\ref{RemnormaVAr}, we have that
\Beq
  \< \calR_r v , w >_{A,r}
  = (v,e_1) (w,e_1) + (A^r v,A^r w)
  = |\Omega| \, (\mean v)(\mean w) + (A^r v,A^r w)
  \non
\Eeq
for every $v,w\in\VA r$.
In particular, if $v\in\Vz r$, then we have $\mean v=0$.
Moreover, $A^r(1)=0$ since $\lambda_1=0$.
Hence,
\Beq
  \< \calR_r v , 1 >_{A,r}
  = (A^r v,A^r(1))
  = 0
  \quad \hbox{for every $w\in\VA r$}.
  \non  
\Eeq
This shows that $\calR_r v\in\Vz{-r}$ for every $v\in\Vz r$.
Now, we fix any $f\in\Vz{-r}$ and prove that the element
$\,v:=\calR_r^{-1}f\,$ of $\,\VA r\,$ belongs to~$\Vz r$.
We have, indeed,
\Beq
  0 = \< f,1 >_{A,r}
  = \< \calR_r v,1 >_{A,r}
  = |\Omega| \, (\mean v)(\mean 1) + (A^r v,A^r(1))
  = |\Omega| \, \mean v.
  \non
\Eeq
This concludes the proof of the first assertion of the statement.
The second one means that, for every $v\in\Vz{2r}$, 
the elements $\calR_r v\in\VA{-r}$ and $A^{2r}v\in H$ coincide
in the sense of the embedding $H\subset\VA{-r}$.
Thus, we fix $v\in\Vz{2r}$ and $w\in\VA r$. 
In both cases
$\lambda_1>0$ and $\lambda_1=0$ 
(in~the latter since $\mean v=0$),
we~have by the definition \eqref{defprodAr} of the inner product that
\Bsist
  && \< \calR_r v,w >_{A,r}
  = (A^r v,A^r w)
  = \somma j1\infty \bigl( \lambda_j^r (v,e_j) \bigr) \bigl( \lambda_j^r (w,e_j) \bigr)
  \non
  \\
  && = \somma j1\infty \bigl( \lambda_j^{2r}(v,e_j) \bigr) (w,e_j)
  = (A^{2r} v,w)
  = \< A^{2r} v,w >_{A,r} \,.
  \non
\Esist
As $w\in\VA r$ is arbitrary, we conclude that $\calR_r v=A^{2r}v$.
\Edim

Due to the above result, it is reasonable to use a proper notation
for the restrictions of $\calR_r$ and $\calR_r^{-1}$ to the subspaces
$\Vz r$ and $\Vz{-r}$, respectively.
We~set 
\Beq
  \Az{2r} := (\calR_r)_{|\Vz r}
  \aand
  \Az{-2r} := (\calR_r^{-1})_{|\Vz{-r}}\,,
  \label{defAz}
\Eeq
where the index $0$ means nothing if $\lambda_1>0$ (since then $\Vz{\pm r}=\VA{\pm r}$),
while it reminds the zero mean value condition in the case $\lambda_1=0$.
We thus have
\Bsist
  && \Az{2r} \in \calL(\Vz r,\Vz{-r}) , \quad
  \Az{-2r} \in \calL(\Vz{-r},\Vz r)
  \aand
  \Az{-2r} = (\Az{2r})^{-1}\,,
  \label{contlinAz}
  \\[1mm]
  && \< \Az{2r} v,w >_{A,r} = (v,w)_{A,r} = (A^r v,A^r w)
  \quad \hbox{for every $v\in\Vz r$ and $w\in\VA r$}\,,
  \qquad
  \label{identityAz}
  \\[1mm]
  && \< f , \Az{-2r} f >_{A,r}
  = \norma{\Az{-2r} f}_{A,r}^2
  = \norma f_{A,-r}^2
  \quad \hbox{for every $f\in\Vz{-r}$}.
  \label{normaAz}
\Esist
Notice that \eqref{normaAz} implies that
\Beq
  \< f' , \Az{-2r} f >_{A,r}
  = \frac 12 \, \frac d {dt} \,\norma f_{A,-r}^2
  \quad \hbox{\aet, \ for every $f\in\H1{\Vz{-r}}$}.
  \label{dtnormaAz}
\Eeq

\Bprop
\label{Ok}
We have
\Beq
  \bigl( A^r \Az{-2r} f , A^r v )
  = \< f,v >_{A,r}
  \quad \hbox{for every $f\in\Vz{-r}$ and $v\in\VA r$}.
  \label{ok}
\Eeq
\Eprop

\Bdim
We first notice that 
$(e_i,e_j)_{A,r}=(\lambda_i^re_i,\lambda_j^re_j)=\lambda_j^{2r}\delta_{ij}$
for $i,j\geq2$,
so that the system $\{\lambda_j^{-r}e_j\}_{j\geq2}$ is orthonormal in~$\VA r$.
It follows that
\Bsist
  && \hbox{the series} \quad
  \textstyle\somma j2\infty c_j e_j
  = \somma j2\infty (\lambda_j^r c_j) (\lambda_j^{-r} e_j)
  \quad \hbox{converges in $\VA r$}
  \non
  \\[1mm]
  && \hbox{if and only if} \quad
  \textstyle\somma j2\infty |\lambda_j^r c_j|^2 < + \infty 
  \quad \hbox{or} \quad
  \textstyle\somma j1\infty |\lambda_j^r c_j|^2 < + \infty \,.
  \non
\Esist
On the other hand, if $v\in\VA r$, we have both
\Beq
  \somma j1\infty |\lambda_j^r (v,e_j)|^2 < + \infty
  \aand
  \somma j1\infty (v,e_j) e_j = v
  \quad \hbox{in $H$}.
  \non
\Eeq
We conclude that
\Beq
  \somma j1\infty (v,e_j) e_j = v
  \quad \hbox{in $\VA r$}
  \quad \hbox{for every $v\in\VA r$} .
  \label{fourierVAr}
\Eeq
In particular, if we set, for convenience,
\Beq
  \jz = 1
  \quad \hbox{if $\lambda_1>0$}
  \aand
  \jz = 2
  \quad \hbox{if $\lambda_1=0$,}
  \label{defjz}
\Eeq
then we have that
\Beq
  \somma j\jz\infty (v,e_j) e_j = v
  \quad \hbox{in $\Vz r$}
  \quad \hbox{for every $v\in\Vz r$} .
  \label{fourierVzr}
\Eeq
Next, notice that $e_j \in \VA{2r} \cap \Vz r = \Vz{2r}$ if $j\geq\jz$, 
whence, by Proposition~\ref{ExtensionR},
\Beq
  \Az{2r} e_j = \calR_r e_j = A^{2r} e_j = \lambda_j^{2r} e_j
  \quad \hbox{for every $j\geq\jz$}.
  \non
\Eeq
Now, take any $f\in\Vz{-r}$ and set $z:=\Az{-2r}f$.
Then $z\in\Vz r$ so that \eqref{fourierVzr} holds for~$z$.
Therefore, since $\Az{2r}\in\calL(\Vz r,\Vz{-r})$, we deduce that
\Beq
  f = \Az{2r} z
  = \somma j\jz\infty (z,e_j) \Az{2r} e_j
  = \somma j\jz\infty \lambda_j^{2r} (z,e_j) e_j
  \quad \hbox{in $\Vz{-r}$},
  \non
\Eeq
whence also
\Beq
  \< f,e_i >_{A,r}
  = \somma j\jz\infty \lambda_j^{2r} (z,e_j) \< e_j , e_i >_{A,r}
  = \somma j\jz\infty \lambda_j^{2r} (z,e_j) ( e_j , e_i )
  = \lambda_i^{2r} (z,e_i)
  \quad \hbox{for every $i\geq\jz$}.
  \non
\Eeq
Hence, the above series expansion becomes
\Beq
  f = \somma j\jz\infty \< f,e_j>_{A,r} \, e_j
  \quad \hbox{in $\Vz{-r}$}
  \quad \hbox{for every $f\in\Vz{-r}$}.
  \label{fourierVzmenor}
\Eeq
At this point, we can easily conclude.
Indeed, on the one side, the formulas \eqref{fourierVzmenor} and 
\eqref{fourierVAr}, 
combined with $\Az{-2r}\in\calL(\Vz{-r},\Vz r)$ and $A^r\in\calL(\VA r,H)$,
ensure that
\Bsist
  && \bigl( A^r \Az{-2r} f , A^r v )
  = \Bigl( A^r \Az{-2r} \textstyle\somma j\jz\infty \< f,e_j >_{A,r} e_j , A^r \textstyle\somma j1\infty (v,e_j) e_j \Bigr)
  \non
  \\
  && = \Bigl( \textstyle\somma j\jz\infty \< f,e_j >_{A,r} \, A^r \Az{-2r} e_j , \textstyle\somma j1\infty (v,e_j) A^r e_j \Bigr)
  \non
  \\
  && = \Bigl( \textstyle\somma j\jz\infty \< f,e_j >_{A,r} \lambda_j^{-r} e_j , \textstyle\somma j1\infty (v,e_j) \lambda_j^r e_j \Bigr)
  = \somma j\jz\infty \< f,e_j >_{A,r} (v,e_j)
  \non
  \\
  && = \< f , \textstyle\somma j\jz\infty (v,e_j) e_j >_{A,r} 
    \quad \pier{\hbox{for every $f\in\Vz{-r}, \ \, v\in V_A^r $}}.
  \non
\Esist
On the other hand, the last expression is equal to $\<f,v>_{A,r}$
in both the cases $\lambda_1>0$ and $\lambda_1=0$, 
since the assumption $f\in\Vz{-r}$ implies that $\<f,1>_{A,r}=0$ \,in the latter.
\Edim

\Bprop
\label{ExtensionA}
The operator $A^{2r}\in\calL(\VA{2r},H)$ can be extended in a unique way
to a continuous linear operator, still termed~$A^{2r}$, from $\VA r$ into~$\Vz{-r}$.
Moreover,
\Beq
  \norma{A^{2r}v}_{A,-r} \leq \norma{A^r v}
  \quad \hbox{for every $v\in\VA r$}.
  \label{stimaAdr}
\Eeq
\Eprop

\Bdim
For $v\in\VA{2r}$ and $w\in\VA r$, we have that
\Bsist
  && \< A^{2r}v , w >_{A,r}
  = (A^{2r}v , w)
  = \somma j1\infty ( \lambda_j^{2r} (v,e_j) ) (w,e_j)
  \non
  \\
  && = \somma j1\infty ( \lambda_j^r (v,e_j) ) ( \lambda_j^r (w,e_j) )
  = (A^r v,A^r w)
  \leq \norma v_{A,r} \norma w_{A,r} \,.
  \non
\Esist
We deduce that
\Beq
  \norma{A^{2r}v}_{A,-r}
  \leq \norma v_{A,r}
  \quad \hbox{for every $v\in\VA{2r}$}.\non
\Eeq
This shows that the mapping\, $\VA{2r}\ni v\mapsto A^{2r}v\in\VA{-r}$\,
is continuous if $\VA{2r}$ is endowed with the topology induced by~$\VA r$.
On the other hand, $\VA{2r}$~is dense in~$\VA r$
(see~\eqref{compembA}).
Thus, the existence of a unique extension
$A^{2r}\in\calL(\VA r,\VA{-r})$ follows,
and we have
\Beq
  \< A^{2r}v , w >_{A,r} 
  = (A^r v , A^r w)
  \quad \hbox{for every $v,w\in\VA r$}.
  \label{extensionA}
\Eeq
We immediately infer that
\Beq
  |\< A^{2r}v , w >_{A,r}|
  \leq \norma{A^r v} \, \norma{A^r w}
  \leq \norma{A^r v} \, \norma w_{A,r}
  \quad \hbox{for every $v,w\in\VA r$},
  \non
\Eeq
whence \eqref{stimaAdr} clearly follows.
Thus, it remains to verify that $A^{2r}v\in\Vz{-r}$ for every $v\in\VA r$ if $\lambda_1=0$
(since there is nothing to prove if $\lambda_1>0$).
For every $v\in\VA r$, we have
\Beq
  \< A^{2r}v , 1 >_{A,r} 
  = (A^r v , A^r(1))
  = 0 ,
  \non
\Eeq
since $\lambda_1=0$ implies that $A^r(1)=0$ by \eqref{hpsimple}.
Hence, \pier{it turns out that} $A^{2r}v\in\Vz{-r}$, as claimed.
\Edim

\Bprop
\label{NormeVA}
For every $f\in\VA{-r}$, we have the representations
\Bsist
  && \norma f_{A,-r}^2
  = \somma j1\infty |\lambda_j^{-r} \<f,e_j>_{A,r}|^2
  \quad \hbox{if $\lambda_1>0$},
  \qquad
  \label{normaVA}
  \\[-5pt]
  && \norma f_{A,-r}^2
  = |\<f,e_1>_{A,r}|^2 + \somma j2\infty |\lambda_j^{-r} \<f,e_j>_{A,r}|^2
  \quad \hbox{if $\lambda_1=0$}.
  \label{normaVAbis}
\Esist
\Eprop

\Bdim
Assume first that $\lambda_1>0$ and set $w:=\calR_r^{-1}f$.
Then the definition \eqref{defnormaAr} yields that
\Beq
  \norma w_{A,r}^2
  = \somma j1\infty |\lambda_j^r (w,e_j)|^2 \,.
  \non
\Eeq
On the other hand, by the definition of the Riesz operator~$\calR_r$, we have that
\Beq
  \< \calR_r w , v >_{A,r}
  = (w,v)_{A,r}
  = (A^r w,A^r v)
  = \somma j1\infty \lambda_j^{2r} (w,e_j) (v,e_j)
  \quad \hbox{for every $v\in\VA r$}.
  \non
\Eeq
In particular, it also holds the identity
\Beq
  \< f,e_i>_{A,r}
  = \lambda_i^{2r} (w,e_i)
  \quad \hbox{for every $i\geq1$}.
  \non
\Eeq
Therefore, by recalling \eqref{defnormaAr}, we deduce that
\Beq
  \norma f_{A,-r}^2
  = \norma w_{A,r}^2
  = \somma j1\infty |\lambda_j^r (w,e_j)|^2
  = \somma j1\infty |\lambda_j^r \, \lambda_j^{-2r} \<f,e_j>_{A,r}|^2
  = \somma j1\infty |\lambda_j^{-r} \<f,e_j>_{A,r}|^2\,,
  \non
\Eeq
that is, \eqref{normaVA} is valid.
If, instead, $\lambda_1=0$, then
the same calculation with $\lambda_1$ replaced by~$1$ yields that
\Beq
  \norma f_{A,-r}^2  
  = |(w,e_1)|^2 + \somma j2\infty |\lambda_j^{-r} \<f,e_j>_{A,r}|^2 .
  \non
\Eeq
On the other hand, we have that
\Beq
  \< f,e_1>_{A,r}
  = \< \calR_r w,e_1>_{A,r}
  = (w,e_1)_{A,r}
  = (w,e_1) (e_1,e_1) + (A^r w,A^r e_1)
  = (w,e_1),
  \non
\Eeq
since $A^r e_1=0$.
Therefore, \eqref{normaVAbis} follows as well.
\Edim

\Bprop
\label{Interpolation}
For every $\eta>0$ and $v\in\VA\eta$,
there holds the interpolation inequality
\Beq
  \norma v \leq \norma v_{A,\eta}^\theta \, \norma v_{A,-r}^{1-\theta},
  \quad \hbox{where} \quad
  \theta = \frac r {r+\eta} \,.
  \label{interpolation}
\Eeq
\Eprop

\Bdim
Set $c_j:=(v,e_j)$ for $j\geq1$, for brevity,
and first assume that $\lambda_1>0$.
Then we have
\Beq
  \norma v_{A,\eta}^2 = \somma j1\infty |\lambda_j^\eta c_j|^2
  \aand
  \norma v_{A,-r}^2 = \somma j1\infty |\lambda_j^{-r} c_j|^2,
  \non
\Eeq
thanks to \eqref{normaVA}.
Therefore, by using the H\"older inequality for infinite sums
and noticing that $(1-\theta)r/\theta=\eta$, we find that
\Bsist
  && \norma v^2
  = \somma j1\infty c_j^2
  = \somma j1\infty \lambda_j^{2(1-\theta)r} c_j^{2\theta} \lambda_j^{-2(1-\theta)r} c_j^{2(1-\theta)}
  \non
  \\
  && \leq \Bigl( \somma j1\infty |\lambda_j^{2(1-\theta)r} c_j^{2\theta}|^{\frac 1\theta} \Bigr)^\theta
    \Bigl( \somma j1\infty |\lambda_j^{-2(1-\theta)r} c_j^{2(1-\theta)}|^{\frac 1{1-\theta}} \Bigr)^{1-\theta}
  \non
  \\
  && = \Bigl( \somma j1\infty |\lambda_j^\eta c_j|^2 \Bigr)^\theta
    \Bigl( \somma j1\infty |\lambda_j^{-r} c_j|^2 \Bigr)^{1-\theta}
  = \norma v_{A,\eta}^{2\theta} \, \norma v_{A,-r}^{2(1-\theta)} \,.
  \non
\Esist
Assume now that $\lambda_1=0$.
Then the same calculation with $\lambda_1$ replaced by $1$ yields that
\Beq
  \norma v^2
  = \somma j1\infty c_j^2
  \leq \Bigl( c_1^2 + \somma j2\infty |\lambda_j^\eta c_j|^2 \Bigr)^\theta
    \Bigl( c_1^2 + \somma j2\infty |\lambda_j^{-r} c_j|^2 \Bigr)^{1-\theta}
  = \norma v_{A,\eta}^{2\theta} \, \norma v_{A,-r}^{2(1-\theta)}.
  \non
\Eeq
Hence, the inequality \eqref{interpolation} holds true in any case.
\Edim

\Brem
\label{InterpolationzeroT}
By simply applying the above result and owing to the H\"older inequality,
we deduce that
\Beq
  \norma v_{\L2H}
  \leq \norma v_{\L2{\VA\eta}}^\theta \, \norma v_{\L2{\VA{-r}}}^{1-\theta}
  \quad \hbox{for every $v\in\L2{\VA\eta}$},
  \label{interpolationzeroT}
\Eeq
with the same $\theta$ as in~\eqref{interpolation}.
\Erem

Now, we introduce some notations concerning interpolating functions.

\Bnot
\label{Interpolants}
Let $N$ be a positive integer and $Z$ be one of the spaces $H$, $\VA r$,~$\VB\sigma$.
We set $h:=T/N$ and $I_n:=((n-1)h,nh)$ for $n=1,\dots,N$.
Given $z=(z^0,z^1,\dots ,z^N)\in Z^{N+1}$,
we define the piecewise constant and piecewise linear interpolants
\Beq
  \overline z_h \in \L\infty Z , \quad
  \underline z_h \in \L\infty Z 
  \aand
  \hat z_h \in \W{1,\infty}Z
  \non
\Eeq
by setting 
\Bsist
  && \hskip -2em
  \overline z_h(t) = z^n
  \aand
  \underline z_h(t) = z^{n-1}
  \quad \hbox{for a.a.\ $t\in I_n$, \ $n=1,\dots,N$},
  \label{pwconstant}
  \\
  && \hskip -2em
  \hat z_h(0) = z_0
  \aand
  \dt\hat z_h(t) = \frac {\pier{z^{n}-z^{n-1}}} h
  \quad \hbox{for a.a.\ $t\in I_n$, \ $n=1,\dots,N$}.
  \qquad
  \label{pwlinear}
\Esist
\Enot

For the reader's convenience,
we summarize the relations between the finite set of values
and the interpolants in the following proposition,
whose proof follows from \sfw\ computation:

\Bprop
\label{Propinterp}
With Notation~\ref{Interpolants}, we have that
\Bsist
  && \norma{\overline z_h}_{\L\infty Z}
  = \max_{n=1,\dots,N} \norma{z^n}_Z \,, \quad
  \norma{\underline z_h}_{\L\infty Z}
   = \max_{n=0,\dots,N-1} \norma{z^n}_Z\,,
  \label{ouLinftyZ}
  \\
  && \norma{\dt\hat z_h}_{\L\infty Z}
  = \max_{\pier{n=0,\dots,N-1}}\norma{(z^{n+1}-z^n)/h}_Z\,,
  \label{dtzLinftyZ}
  \\
  \separa
  && \norma{\overline z_h}_{\L2Z}^2
  = h \somma n1N \norma{z^n}_Z^2 \,, \quad
  \norma{\underline z_h}_{\L2Z}^2
  = h \somma n0{N-1} \norma{z^n}_Z^2 \,,
  \label{ouLdueZ}
  \\
  \separa
  && \norma{\dt\hat z_h}_{\L2Z}^2
  = h \somma n0{N-1} \norma{(z^{n+1}-z^n)/h}_Z^2\,, 
  \label{dtzLdueZ}
  \\
  && \norma{\hat z_h}_{\L\infty Z}
  = \max_{n=1,\dots,N} \max\{\norma{z^{n-1}}_Z,\norma{z^n}_Z\}
  = \max\{\norma{z_0}_Z,\norma{\overline z_h}_{\L\infty Z}\}\,,
  \qquad\qquad
  \label{hzLinftyZ}
  \\
  && \norma{\hat z_h}_{\L2Z}^2
  \leq h \somma n1N \bigl( \norma{z^{n-1}}_Z^2 + \norma{z^n}_Z^2 \bigr)
  \leq h \norma{z_0}_Z^2
  + 2 \norma{\overline z_h}_{\L2Z}^2 \,.
  \label{hzLdueZ}
\Esist
Moreover, it holds that
\Bsist
  && \norma{\overline z_h-\hat z_h}_{\L\infty Z}
  = \max_{n=0,\dots,N-1} \norma{z^{n+1}-{z^n}}_Z
  = h \, \norma{\dt\hat z_h}_{\L\infty Z}\,,
  \qquad
  \label{diffLinfty}
  \\
  && \norma{\overline z_h-\hat z_h}_{\L2Z}^2
  = \frac {{h}} 3 \somma n0{N-1} \norma{z^{n+1}-z^n}_Z^2
  = \frac {h^2} 3 \, \norma{\dt\hat z_h}_{\L2Z}^2\,,
  \label{diffLdue}
\Esist
and similar identities for the difference $\underline z_h-\hat z_h$.
As a consequence, we have the inequalities
\Bsist
  && \norma{\overline z_h-\underline z_h}_{\L\infty Z}
  \leq 2h \, \norma{\dt\hat z_h}_{\L\infty Z}\,,
  \qquad
  \label{diffbisLinfty}
  \\
  && \norma{\overline z_h-\underline z_h}_{\L2Z}^2
  \leq \frac {\pier{4} h^2} 3 \, \norma{\dt\hat z_h}_{\L2Z}^2 \,.
  \label{diffbisLdue}
\Esist
Finally, we have that
\Bsist
  && h \somma n0{N-1} \norma{(z^{n+1}-z^n)/h}_Z^2 
  \leq \norma{\dt z}_{\L2Z}^2
  \non
  \\
  && \quad \hbox{if $z\in\H1Z$\aand $z^n=z(nh)$ for $n=0,\dots,N$}.
  \label{interpH1Z}
\Esist
\Eprop

\vspace{3mm}
Throughout the paper, we make use of
the elementary identity and inequalities
\Bsist
  \hskip-1cm&& a (a-b)
  = \frac 12 \, a^2
  + \frac 12 \, (a-b)^2
  - \frac 12 \, b^2
  \geq \frac 12 \, a^2
  - \frac 12 \, b^2
  \quad \hbox{for every $a,b\in\erre$},
  \label{elementare}
  \\
  \hskip-1cm&& ab\leq \delta a^2 + \frac 1 {4\delta}\,b^2
  \quad \hbox{for every $a,b\in\erre$ and $\delta>0$},
  \label{young}
\Esist
\Accorpa\Elementari elementare young
and quote \eqref{young} as the Young inequality.
We also take advantage of the summation by parts formula
\Beq
  \somma n0{k-1} a_{n+1} (b_{n+1} - b_n)
  = a_k b_k - a_1 b_0
  - \somma n1{k-1} (a_{n+1} - a_n) b_n\,,
  \label{byparts}
\Eeq
which is valid for arbitrary real numbers $a_1,\dots,a_k$ and $b_0,\dots,b_k$.
We also account for the discrete Gronwall lemma in the following form
(see, e.g., \cite[Prop.~2.2.1]{Jerome}):
for nonnegative real numbers $M$ and $a_n,b_n$, $n=0,\dots,N$, 
\Bsist
  a_k \leq M + \somma n0{k-1} b_n a_n 
  \quad \hbox{for $k=0,\dots,N$}
  \qquad \hbox{implies}
  \non
  \\
  a_k \leq M \exp \Bigl( \somma n0{k-1} b_n \Bigr)
  \quad \hbox{for $k=0,\dots,N$}.
  \label{gronwall}
\Esist
In \accorpa{byparts}{gronwall} it is understood that
a sum vanishes if the corresponding set of indices is empty.

Finally, we state a general rule that we follow throughout the paper 
as far as the constants are concerned.
We always use a small-case italic $c$ without subscripts
for~different constants that may only depend on the final time~$T$,
the operators $A^r$ and~$B^\sigma$,  
the shape of the nonlinearities $\beta$ and~$\pi$,
and the properties of the data involved in the statements at hand.
Thus, the values of such constants do not depend on~$\tau$,
nor on the regularization parameter $\lambda$ or the time step $h$ we 
introduce in Section~\ref{APPROX},
and it is clear that they might change from line to line 
and even in the same formula or chain of inequalities. 
In contrast, we use different symbols (e.g., capital letters like $M_0$ in~\eqref{hpdatireg})
for precise values of constants we want to refer~to.


\section{Continuous dependence and uniqueness}
\label{UNIQUENESS}
\setcounter{equation}{0}

This section is devoted to the proof of the \jupiter{continuous dependence and uniqueness result stated in Theorem~\ref{Wellposedness}.
We consider only the case of the same initial datum, for simplicity;
however,} the case of different initial data sketched in Remark~\ref{Remcontdep}
could be \previousgianni{treated} in the same way with only minor changes.

We pick two data $u_i$, $i=1,2$, \jupiter{and corresponding solutions} $(y_i,\mu_i)$, 
and set for convenience
$u:=u_1-u_2$, $y:=y_1-y_2$ and $\mu:=\mu_1-\mu_2$.
Now, we write equation \eqref{prima} at the time $s$ for these solutions and take the difference.
Then, we test it by $v=\Az{-2r}y(s)$ by observing that $y(s)\in\Vz{-r}$
since $y\in\L2H$ and $\mean y(s)=0$ if $\lambda_1=0$ by the conservation property~\eqref{conservation},
so that $v$ is a well-defined element of~$\VA r$.
Moreover, $\Az{-2r}y\in\L\infty{\VA r}$, since $y\in\L\infty{\VA{-r}}$ by~\eqref{regy}.
Integrating over  $(0,t)$ with respect to~$s$, where $t\in(0,T)$ is arbitrary,
we obtain the identity
\Beq
  \iot \< \dt y(s) , \Az{-2r} y(s) >_{A,r} \, ds
  + \iot \bigl( A^r \mu(s) , A^r \Az{-2r}y(s) \bigr) \, ds
  = 0 \,.
  \non
\Eeq
Now, we apply~\eqref{dtnormaAz} and \eqref{ok}, 
noting that $\mu\in\L2{\VA r}$.
Thus, the above identity becomes
\Beq
  \frac 12 \, \norma{y(t)}_{A,-r}^2
  + \iot (y(s) , \mu(s)) \, ds 
  = 0
  \label{testprima}
\Eeq
where the duality product of \eqref{ok} has been replaced
by the inner product here, since both $y$ and $\mu$ are $H$-valued.
At the same time, we write \eqref{seconda} for $u_i$ and $(y_i,\mu_i)$, $i=1,2$,
test them by $y_2$ and~$y_1$, respectively, 
add the resulting inequalities to each other, and integrate over $(0,t)$ as before.
Then, the terms involving $\Beta$ cancel out, and we obtain (after rearranging) that
\Bsist
  && \frac \tau 2 \, \norma{y(t)}^2 
  + \iot \norma{B^\sigma y(s)}^2\,ds
  - \iot \bigl( \mu(s) , y(s) \bigr) \, ds
  \non
  \\
  && \leq \iot \bigl( u(s),y(s) \bigr) \, ds\,
  - \iot \bigl( \pi(y_1(s)) - \pi(y_2(s)) , y(s) \bigr) \, ds \,.
  \non
\Esist
By adding this to \eqref{testprima}, 
and accounting for the \Lip\ continuity of~$\pi$ and the Schwarz and Young inequalities,
we deduce that (with $\Lpi'$ given by~\eqref{defbetapi})
\Beq
  \frac 12 \, \norma{y(t)}_{A,-r}^2
  + \frac \tau 2 \, \norma{y(t)}^2 
  + \iot \norma{B^\sigma y(s)}^2\,ds
  \leq \frac 14 \iot \norma{u(s)}^2 \, ds
  + \Lpi' \iot \norma{y(s)}^2 \, ds.
  \label{forcontdep}
\Eeq
At this point, we recall the compacness inequality~\eqref{compineq}.
Thus, we have that
\Beq
  \Lpi' \iot \norma{y(s)}^2 \, ds
  \leq \frac 12 \iot \norma{B^\sigma y(s)}^2 \, ds
  + c \iot \norma{y(s)}_{A,-r}^2 \, ds .
  \non
\Eeq
By combining this with \eqref{forcontdep} and applying the Gronwall lemma,
we conclude that the desired estimate \eqref{contdep} holds true
with a constant $K_2$ as in the statement. \jupiter{In particular, it turns out that
$y_1=y_2$ if $u_1=u_2$. In this case, testing the difference of \eqref{prima}, 
written for the solutions, by $v= \mu : = \mu_1-\mu_2$, we deduce that $A^r \mu (t) = 0 $
\aat. If $\lambda_1 >0$, this implies that $ \mu (t) =0 $ \aat. This concludes the proof 
of the theorem.}

\Brem
\label{Jupiter} 
\jupiter{Let us point out at this place that the uniqueness property for the
second component $\mu$ of the solutions $(y,\mu) $ to \eqref{regy}--\eqref{cauchy} can also 
hold true in other situations if $\lambda_1=0$. For instance, suppose that $\Beta \in C^1 (\erre)$ and $y\in L^\infty (Q)$. Then, in view of \eqref{hpVB}, we may, for a constant $\delta>0$, insert $v=y-\delta$ 
in the variational inequality \eqref{seconda} for $\mu_1$ and $v=y+\delta$ 
in \eqref{seconda} for $\mu_2$, where $(y,\mu_1)$ and $(y,\mu_2)$ are two solutions. After addition of the two resulting inequalities, we then obtain \aet\ that 
$$
2\iO \Beta ( y) \leq \delta (\mu_1 - \mu_2, 1) + \iO \Beta ( y-\delta) + \iO \Beta (y+\delta), 
$$ 
where the two integrals on the \rhs\ are finite. Division by $\delta >0$ then yields that
$$
\iO \frac{\Beta ( y) - \Beta ( y-\delta)}\delta  + \iO \frac{\Beta ( y) - \Beta 
( y+\delta)}\delta
\leq (\mu_1 - \mu_2, 1) . 
$$
Taking the limit as $\delta \searrow 0$, we conclude from the Lebesgue theorem of dominated convergence that 
$$0= \iO \beta (y) - \iO \beta (y) \leq (\mu_1 - \mu_2, 1). $$
Interchanging the roles of $\mu_1$  and $\mu_2$, we then infer that $\mean \mu_1 = \mean \mu_2$ \aet. Since $A^r(\mu_1 - \mu_2) = 0,$ this implies that $\mu_1 = \mu_2$.  
}
\Erem


\section{Approximation}
\label{APPROX}
\setcounter{equation}{0}

In this section we deal with an approximation of problem \Pbl\ 
and solve it by a time discretization procedure.
We first introduce the \pier{Moreau--Yosida} regularizations 
$\Betal$ and $\betal$ of $\Beta$ of $\beta$ at the level $\lambda>0$
(see, e.g., \cite[p.~28 and p.~39]{Brezis}).
By~accounting for assumptions~\accorpa{hpBeta}{hpcoerc}, we~have
\Bsist
  && \Betal(s) = \int_0^s \betal(s') \, ds'
  \aand
  0 \leq \Betal(s) \leq \Beta(s)
  \quad \hbox{for every $s\in\erre$}.
  \label{defBetal}
  \\[1mm] 
  && \Betal(s) + \Pi(s)
  \geq \alpha \, s^2 - C
  \non
  \\
  &&  \hbox{for some constants $\alpha,C>0$, every $s\in\erre$ and $\lambda>0$ small enough}.
  \label{coerclambda}
\Esist
Moreover, we recall that $\betal$ is \Lip\ continuous, so that $\Betal$ 
grows \juerg{at most} quadratically,
and that the following properties hold true:
\Bsist
  && \Beta_{\lambda'}(s) \geq \Beta_{\lambda''}(s)
  \quad \hbox{if $\lambda'\leq\lambda''$}
  \aand
  \lim_{\lambda\searrow0} \Betal(s) = \Beta(s)
  \quad \hbox{for every $s\in\erre$},
  \label{propBetal}
  \\
  && |\betal(s)| \leq |\betaz(s)|
  \quad \hbox{for every $s\in D(\beta)$}.
  \label{propbetal}
\Esist
By replacing $\Beta$ in \eqref{seconda} by~$\Betal$,
we obtain the following system:
\Bsist
  && \< \dt\yl(t) , v >_{A,r}
  + ( A^r \mul(t) , A^r v )
  = 0
  \quad \hbox{for every $v\in\VA r$ and \aat},
  \qquad
  \label{primal}
  \\[1mm]
  && \previousgianni{ \bigl( \tau\dt\yl(t) , \yl(t) - v \bigr) }
  + \bigl( B^\sigma\yl(t) , B^\sigma(\yl(t)-v) \bigr)
  \non
  \\
  && \quad {}
  + \iO \Betal(\yl(t))
  + \bigl( \pi(\yl(t)) - u(t) ,  \yl(t)-v \bigr)
  \non
  \\
  && \leq \bigl( \mul(t) ,  \yl(t)-v \bigr)
  + \iO \Betal(v)
  \quad \hbox{for every $v\in\VB\sigma$ and \aat},
  \label{secondal}
  \\[1mm]
  && \yl(0) = \previousgianni\yz . 
  \label{cauchyl}
\Esist
\Accorpa\Pbll primal cauchyl
We stress that \eqref{secondal} is equivalent to both 
the time-integrated variational inequality
\Bsist
  && \previousgianni{\ioT \bigl(\tau  \dt\yl(t) , \yl(t) - v(t) \bigr) \, dt }
  + \ioT \bigl( B^\sigma\yl(t) , B^\sigma(\yl(t)-{v(t)}) \bigr) \, dt
  \non
  \\
  && \quad {}
  + \intQ \Betal(\yl)
  + \ioT \bigl( \pi(\yl(t)) - u(t) ,  \yl(t)-{v(t)} \bigr) \, dt
  \non
  \\
  && \leq \ioT \bigl( \mul(t) ,  \yl(t)-{v(t)} \bigr) \, dt
  + \intQ \Betal(v)
  \quad \hbox{for every $v\in\L2{\VB\sigma}$},
  \label{intsecondal}
\Esist
and \jupiter{the pointwise variational equation}
\Bsist
  && \jupiter{( \tau\dt\yl  (t) , v)} + \bigl( B^\sigma\yl(t) , B^\sigma v \bigr)
  + \bigl( \betal(\yl(t)) + \pi(\yl(t)) - u(t) ,  v \bigr)
  = \bigl( \mul(t) ,  v \bigr)
  \non
  \\[1mm]
  &&\quad \hbox{for every $v\in\VB\sigma$ and \aat}.
  \label{eqsecondal}
\Esist
\jupiter{The derivation of \eqref{eqsecondal} from \eqref{secondal} uses the fact that 
$\Betal $ is differentiable and $\betal$ is its Lipschitz continuous derivative. From \eqref{eqsecondal} we can infer that $B^{2\sigma } y_\lambda \in\L2H$, so that
$$
\mul = \tau \dt\yl + B^{2\sigma} \yl+ \betal(\yl) + \pi(\yl) - u 
$$
is uniquely determined. We therefore conclude from Theorem~\ref{Wellposedness} the following result.}

\Bthm
\label{Wellposednessl}
Under the assumptions of Theorem~\ref{Wellposedness},
problem \Pbll\ has a unique solution satisfying \accorpa{regy}{regmu}.
\Ethm

Uniqueness follows from Theorem~\ref{Wellposedness},
since $\betal$ and $\Betal$ satisfy the properties 
we have postulated for $\beta$ and~$\Beta$.
So, we just have to prove the existence of a solution,
and the remainder of the section is devoted to this proof.
To this end, we solve a proper discrete problem and take the limits of the interpolants
as the time step tends to zero.

\step
The discrete problem

We fix an integer $N>1$ and set $h:=T/N$.
Then, the discrete problem consists in finding
two $(N+1)$-tuples\, $(y^0,\dots,y^N)$ and $(\mu^0,\dots,\mu^N)$
satisfying
\Bsist
  y^0 = \yz \,, \quad
  \mu^0 = 0 , \quad
  (\yu,\dots,y^N) \in (\VB{2\sigma})^N
  \aand
  (\muu,\dots,\mu^N) \in (\VA{2r})^N
  \label{regdiscr}
\Esist
and solving
\begin{align}
  & \dhyn + \munp + A^{2r} \munp
  \,=\, \mun,
  \label{primad}
  \\[1mm]
  & \tau \, \dhyn
  + (\Lpi' I + B^{2\sigma} + \betal + \pi)(\ynp)
   \,=\, \Lpi' \yn + \munp + \unp,
  \label{secondad}
\end{align}
for $n=0,1,\dots,N-1$, where $I:H\to H$ is the identity, $\Lpi'$~is given by 
\eqref{defbetapi}, and
\Beq
  \un := u(nh)
  \quad \hbox{for $n=0,1,\dots,N$}.
  \label{defun}
\Eeq
This problem can be solved inductively for $n=0,\dots,N-1$ in the following way:
let $(\yn,\mun)$ be given in $\VB\sigma\times\VA{2r}$.
We first rewrite the above equations in the form
\begin{align}
  & h (I+A^{2r}) \munp + \ynp
  \,=\, \yn + h \mun, 
  \label{primadbis}
  \\[1mm]
  & ((\Lpi'+(\tau/h))I + B^{2\sigma} + \betal + \pi)(\ynp)
   \,=\, (\Lpi'+(\tau/h)) \yn + \munp + \unp .
  \label{secondadbis}
\end{align}
Next, we observe that the operator $\,\calA_\lambda:=\Lpi I+\betal+\pi:H\to H\,$
is monotone and continuous.
On the other hand, the unbounded operator $\,B^{2\sigma}\,$~is monotone in $H$,
and $I+ B^{2\sigma}$ is surjective, whence it follows that
 $B^{2\sigma}$ is maximal monotone.
Therefore, the sum $\,\calA_\lambda+B^{2\sigma}\,$ is also maximal monotone
(see, e.g., \cite[Cor.~2.1 p.~35]{Barbu}).
It follows that $(1+(\tau/h))I+\calA_\lambda+B^{2\sigma}$, 
i.e., the operator that acts on $\ynp$ in~\eqref{secondad},
is surjective and one-to-one from $\VB{2\sigma}$ onto~$H$.
Therefore, \eqref{secondad}~can be rewritten in the equivalent form
\Beq
  \ynp
  = (L_hI + B^{2\sigma} + \betal + \pi)^{-1} \bigl(L_h \yn + \munp + \unp),
  \label{secondadter}
\Eeq
where, for brevity, we have set\, $L_h:=\Lpi'+(\tau/h)$.
\pier{By accounting for \eqref{primadbis}, we conclude that problem \accorpa{primad}{secondad}
is equivalent to the system obtained by coupling \eqref{secondadter} with
the equation
\Beq
  h (I+A^{2r}) \munp
  + (L_hI + B^{2\sigma} + \betal + \pi)^{-1} \bigl(L_h \yn + \munp + \unp)
  \,=\, \yn + h \mun .
  \label{secondadquater}
\Eeq
Arguing} as before, we see that the operator acting on $\munp$ 
on the \lhs\ of \eqref{secondadquater} is surjective and one-to-one from $\VA{2r}$ onto~$H$,
so that the equation can be uniquely solved for~$\munp$ in~$\VA{2r}$.
Inserting the solution in \eqref{secondadter}, 
we directly find that $\ynp\in\VB{2\sigma}$.

\medskip

Once the discrete problem is solved, we can start estimating.
According to the general rule stated at the end of Section~\ref{AUX},
the (possibly different) values of the constants termed $c$ 
are independent of the three parameters $h$, $\lambda$ and~$\tau$.
Moreover, we also express the bounds we find in terms of the interpolants.
According to Notation~\ref{Interpolants}, and recalling that $y^0=\yz\in\VB\sigma$ and that $\mu^0=0$
(see \eqref{hpyz} and~\eqref{regdiscr}),
we remark at once that the discrete problem also reads
\begin{align}
  & \yh \in \W{1,\infty}{\VB\sigma} , \quad
  \underyh \in \L\infty{\VB\sigma}
  \aand
  \overyh \in \L\infty{\VB{2\sigma}},
  \qquad
  \label{regyh}
  \\[1mm]
	  & \undermuh \,, \overmuh \in \L\infty{\VA{2r}},
  \label{regmuh}
  \\[1mm]
  & \dt\yh
  + \overmuh
  + A^{2r} \overmuh
  \,=\, \undermuh
  \quad \aet,
  \label{primah}
  \\[1mm]
  & \tau \, \dt\yh
  + (\Lpi' I + B^{2\sigma} + \betal + \pi)(\overyh)
  \, = \,\Lpi' \underyh + \overmuh + \overuh
  \quad \aet ,
  \label{secondah}
  \\[1mm]
  & \yh(0) = \yz \,.
  \label{cauchyh}
\end{align}

\step
First a priori estimate

We test \eqref{primad} and~\eqref{secondad} 
(by~taking the scalar product in~$H$)
by~$h\munp$ and $\ynp-\yn$, respectively,
and add the resulting identities.
Noting an obvious cancellation, we obtain the equation
\begin{align}
  & h (\munp-\mun,\munp)
  + h (A^{2r} \munp,\munp)
    + \frac \tau h \, \norma{\ynp-\yn}^2
	\non\\[1mm]	
  &+ (B^{2\sigma}\ynp,\ynp-\yn)
    + \bigl( (\Lpi' I + \betal + \pi) (\ynp),\ynp-\yn \bigr)
  \non
  \\[1mm]
  & =\, \Lpi' (\yn,\ynp-\yn)
  + (\unp,\ynp-\yn).
  \non
\end{align}
Now, we observe that the function $\,r\mapsto\frac{\Lpi'}2\,r^2+\previousgianni\Betal(r)+\previousgianni\Pi(r)\,$
is convex on~$\erre$, since $\Betal$ is convex and $\,|\pi'|\leq\Lpi$.
Thus, we have that
\Bsist
  && \bigl( (\Lpi' I + \previousgianni\betal + \previousgianni\pi)(\ynp) , \ynp-\yn \bigr)
  \non
  \\  
  && \geq \frac {\Lpi'} 2 \, \norma\ynp^2 + \iO \bigl( \Betal(\ynp) + \Pi(\ynp) \bigr)
  - \frac {\Lpi'} 2 \, \norma\yn^2 - \iO \bigl( \Betal(\yn) + \Pi(\yn) \bigr).
  \non
\Esist
\gianni{By using this inequality and formulas \accorpa{propA}{propB},
and applying the identity \eqref{elementare}
in two terms on the \lhs\ and in the first one on the \rhs},
we deduce~that
\Bsist
  && \frac h2 \, \norma\munp^2
  + \frac h2 \, \norma{\munp-\mun}^2
  - \frac h2 \, \norma\mun^2
  + h \norma{A^r\munp}^2
  \non
  \\
  && \quad {}
  + \frac \tau h \, \norma{\ynp-\yn}^2
  + \frac 12 \, \norma{B^\sigma \ynp}^2
  + \frac 12 \, \norma{B^\sigma(\ynp-\yn)}^2
  - \frac 12 \, \norma{B^\sigma\yn}^2
  \non
  \\
  && \quad {}
  + \frac {\Lpi'} 2 \, \norma\ynp^2 + \iO \bigl( \Betal(\ynp) + \Pi(\ynp) \bigr)
  - \frac {\Lpi'} 2 \, \norma\yn^2 - \iO \bigl( \Betal(\yn) + \Pi(\yn) \bigr)
  \non
  \\
  \separa
  && \leq \frac {\Lpi'} 2 \, \norma\ynp^2
  - \frac {\Lpi'} 2 \, \norma\yn^2
  \gianni{{} - \frac {\Lpi'} 2 \, \norma{\ynp-\yn}^2}
  + (\unp , \ynp-\yn).
  \non
\Esist
Then, we first \gianni{rearrange} and then sum up for $n=0,\dots,k-1$ with $k\leq N$,
employing summation by parts (see~\eqref{byparts}) in the last term.
\pier{Using \eqref{propBetal}, we} then arrive at the inequality
\begin{align}
  & \frac h2 \, \norma\muk^2
  + \somma n0{k-1} \frac h2 \, \norma{\munp-\mun}^2
  + \somma n0{k-1} h \norma{A^r\munp}^2
  \non
  \\
  & \quad {}
  + \tau \somma n0{k-1} h \Norma{\dhyn}^2
  + \frac 12 \, \norma{B^\sigma \yk}^2
  + \somma n0{k-1} \frac 12 \, \norma{B^\sigma(\ynp-\yn)}^2
  \non
  \\
  & \quad {}
  + \iO \bigl( \Betal(\yk) + \Pi(\yk) \bigr)
  \gianni{{} + \frac {\Lpi'} 2 \somma n0{k-1} \norma{\ynp-\yn}^2}
  \non
  \\
  \separa
  &  \leq  \revis{\frac 12 \, \norma{B^\sigma \yz }^2 
  + \iO \bigl( \Beta(\yz) + \Pi(\yz) \bigr)}
 \non
 \\
  & \quad 
  + (\uk,\yk) - (u^1,\yz) - \somma n1{k-1} (\unp-\un,\yn).
  \label{perprimastima}
\end{align}
Now, we observe that \eqref{coerclambda} implies that
\Beq
  \iO \bigl( \Betal(\yk) + \Pi(\yk) \bigr)
  \geq \frac 12 \iO \bigl( \Betal(\yk) + \Pi(\yk) \bigr)
  + \frac \alpha 2 \, \norma\yk^2 - c\,,
  \non
\Eeq
for \previousgianni{sufficiently} small $\lambda>0$.
In particular, the above integral is bounded from below.
We treat the \rhs\ of \eqref{perprimastima} 
by using the Young and Schwarz inequalities for finite sums, 
as well as~\eqref{interpH1Z}. \revis{Then, for the last three terms of \eqref{perprimastima} we obtain}
\Bsist
  && (\uk,\yk) - (u^1,\yz) - \somma n1{k-1} (\unp-\un,\yn)
  \non
  \\
  && \leq \frac \alpha 4 \, \norma\yk^2 
  + c \, \norma\uk^2
  + \norma\yz^2
  + \norma{u^1}^2
  + \somma n1{k-1} h \Norma\dhun^2
  + \somma n1{k-1} h \norma\yn^2 
  \non
  \\
  && \leq \frac \alpha 4 \, \norma\yk^2
  + \norma\yz^2
  + c \, \norma u_{\L\infty H}^2
  + \norma{\dt u}_{\L2H}^2
  + \somma n1{k-1} h \norma\yn^2 .
  \non
\Esist
By combining the last two estimates with \eqref{perprimastima} and~\eqref{hpyz}
\gianni{and recalling that $\Lpi'\geq1$,}
we infer that
\Bsist
  && \frac h2 \, \norma\muk^2
  + \somma n0{k-1} \frac h2 \, \norma{\munp-\mun}^2
  + \somma n0{k-1} h \norma{A^r\munp}^2
  + \tau \somma n0{k-1} h \Norma{\dhyn}^2
  \non
  \\
  && \quad {}
  + \frac 12 \, \norma{B^\sigma \yk}^2
  + \frac \alpha 4 \, \norma\yk^2
    + \frac 12 \iO \bigl( \Betal(\yk) + \Pi(\yk) \bigr)
  \non
  \\
  && \quad 
  \pier{{}+ \somma n0{k-1} \frac 12 \, \norma{B^\sigma(\ynp-\yn)}^2}
  \gianni{{} + \frac 12 \somma n0{k-1} \norma{\ynp-\yn}^2}
  \non
  \\
  \separa
  && \leq \somma n1{k-1} h \norma\yn^2
  + c \,.
  \non
\Esist
Since this holds for $k=0,\dots,N$,  and as the last integral on the \lhs\ is 
bounded from below,
we~can apply the discrete Gronwall lemma \eqref{gronwall} and conclude that
\Bsist
  && h \, \norma\muk^2
  + \somma n0{k-1} \frac h2 \, \norma{\munp-\mun}^2
  + \somma n0{k-1} h \norma{A^r\munp}^2
  + \tau \somma n0{k-1} h \Norma{\dhyn}^2
  \non
  \\
  && \quad {}
  + \norma\yk_{B,\sigma}^2
  + \iO \bigl( \Betal(\yk) + \Pi(\yk) \bigr)
  \pier{{}+ \somma n0{k-1} \norma{B^\sigma(\ynp-\yn)}^2}
  \gianni{{} + \somma n0{k-1} \norma{\ynp-\yn}^2}
  \non
  \\
  && \leq c 
  \quad\  \hbox{for $k=0,\dots,N$}.
  \label{primastimad}
\Esist
In terms of the interpolants, by neglecting the first contribution
and recalling that $\mu^0=0$,
we have on account of Proposition~\ref{Propinterp} that
\Bsist
  && \norma{\overmuh-\undermuh}_{\L2H}
  + \norma{A^r \overmuh}_{\L2H}
  + \norma{A^r \undermuh}_{\L2H}
  \non
  \\[2pt] 
  && \quad {}
  + \norma\underyh_{\L\infty{\VB\sigma}}
  + \norma\overyh_{\L\infty{\VB\sigma}}
  + \norma\yh_{\L\infty{\VB\sigma}}
  \non
  \\
  && \quad {}
    + \tau^{1/2} \norma{\dt\yh}_{\L2H}
    + \norma{\Betal(\overyh)+\Pi(\overyh)}_{\L\infty\Luno}
  \non
  \\
  && \quad {}
  \pier{{}+ h^{-1/2} \norma{B^\sigma(\overyh-\underyh)}_{\L2H} 
  + h^{-1/2} \norma{\overyh-\underyh}_{\L2H}}
  \leq c \,.
  \label{primastima}
\Esist
\pier{Due to \eqref{hpPi}, we easily infer that} $\,\,\norma{\Pi(\overyh)}_{\L\infty\Luno}\leq c\,(\norma\overyh_{\L\infty H}^2+1)\leq c$, \pier{whence}
we deduce that
\Beq
  \norma{\Betal(\overyh)}_{\L\infty\Luno} \leq c \,.
  \label{daprimastima}
\Eeq

\step
Second a priori estimate

By recalling \eqref{primah} and applying Proposition~\ref{ExtensionA}, we immediately 
obtain
\Bsist
  && \norma{\dt\yh}_{\L2{\VA{-r}}}
  \,\leq\, \norma{\undermuh-\overmuh}_{\L2{\VA{-r}}}
  + \norma{A^{2r}\overmuh}_{\L2{\VA{-r}}}
  \non
  \\[1mm]
  && \leq c \, \norma{\undermuh-\overmuh}_{\L2H}
  + c \, \norma{A^r\overmuh}_{\L2H} \,.
  \non
\Esist
Hence, \eqref{primastima} implies that
\Beq
  \norma{\dt\yh}_{\L2{\VA{-r}}} \leq c \,.
  \label{secondastima}
\Eeq

\step
Consequence

By combining \eqref{primastima} and \eqref{secondastima} with the application of 
\eqref{diffLdue} and its analogue to $\,\overyh$, $\underyh$ and~$\yh$,
we deduce that
\Beq
  \norma{\overyh-\yh}_{\L2{\VA{-r}}}
  + \norma{\underyh-\yh}_{\L2{\VA{-r}}}
  \leq c \, h .
  \label{diffinterpol}
\Eeq

\step 
Third a priori estimate

We want to improve the estimate for $A^r\overmuh$ given by~\eqref{primastima}
and show that
\Beq
  \norma\overmuh_{\L2{\VA r}} 
  + \norma\undermuh_{\L2{\VA r}} \leq c \,.
  \label{terzastima}
\Eeq
By recalling \eqref{hpsimple} and~\eqref{defnormaAr},
we see that there is nothing to prove if $\lambda_1>0$.
On the contrary, if $\lambda_1=0$, we have to estimate the mean value of~$\overmuh$.
Thus, we assume $\lambda_1=0$ and first derive an estimate of~$\betal(\overyh)$.
To this end, \gianni{we recall that $\mz\in\VB\sigma$ by~\eqref{hpVB} 
and that we have postulated the interior assumption~\eqref{hpmz}.
Then, we}
test \eqref{secondad} by~$\,\ynp-\mz\,$ (see~\eqref{hpmz})
and use the inequality
\Beq
  \betal(s) (s-\mz)
  \geq \delta_0 |\betal(s)| - C_0,
  \label{trickMZ}
\Eeq
which holds for some $C_0>0$ and every $s\in\erre$ and $\lambda\in(0,1)$,
where $\delta_0$ is the same as in~\eqref{defdeltaz}
(cf.\ \cite[Appendix, Prop. A.1]{MiZe}; see also \cite[p.~908]{GiMiSchi} for a detailed proof).
By partially using \eqref{primastimad} as well, we have
\begin{align}
  & \iO \bigl( \delta_0 |\betal(\ynp)| - C_0 \bigr)
  \leq \iO \betal(\ynp) (\ynp-\mz)
  \non
  \\
  & = - \tau \Bigl( \dhyn , \ynp-\mz \Bigr)
  - \Lpi' (\ynp-\yn,\ynp-\mz)
  \non
  \\
  & \quad {}
  - \bigl( B^{2\sigma} \ynp , \ynp-\mz \bigr)
  - \bigl( \pi(\ynp) , \ynp-\mz \bigr)
  + (\munp+\unp,\ynp-\mz)
  \non
  \\
  & \leq \pier{c}\,\tau \, \Norma\dhyn \, ( \norma\ynp + 1)
  + c (\norma\ynp^2 + \norma\yn^2 \pier{{}+1})
  + \pier{c}\,\norma\unp \, ( \norma\ynp + 1)
  \non
  \\
  & \quad {}
  + \bigl| \bigl( B^{2\sigma} \ynp , \ynp-\mz \bigr) \bigr|
  + (\munp,\ynp-\mz)
  \non
  \\
  & \leq c \, \tau \, \Norma\dhyn
  + c \, \norma\unp
  + c
  \non
  \\
  & \quad {}
  + \bigl| \bigl( B^{2\sigma} \ynp , \ynp-\mz \bigr) \bigr|
  + (\munp,\ynp-\mz)\,.  
  \label{perterzastima}
\end{align}
We now estimate the last two terms.
For the first one, we owe to assumption \eqref{hpVB} \gianni{just mentioned} and property~\eqref{propB}.
By \pier{recalling} \eqref{primastimad} once more, we see that
\Beq
  \bigl| \bigl( B^{2\sigma} \ynp , \ynp-\mz \bigr) \bigr|
  = \bigl| \bigl( B^\sigma \ynp , B^\sigma\ynp - B^\sigma\mz \bigr) \bigr|
  \leq c \,.
  \non
\Eeq
We deal with the other term by first observing that \eqref{primad} 
and the assumption $\lambda_1=0$ in \eqref{hpsimple} imply that
\Beq
  \mean (\ynp + h \munp) - \mean (\yn + h \mun)
  = \pier{{}-\frac h{|\Omega|}} \, (A^r \munp , A^r(1) ) = 0
  \quad \hbox{},
  \non
\Eeq
\pier{for $n=0,\dots,N-1$,}
whence $\,\mean(\ynp+h\munp)=\mz\,$ for every~$\,n$, since $\mu^0=0$ (see~\eqref{regdiscr}).
Hence, by taking advantage of the Poincar\'e inequality~\eqref{poincare} and~\eqref{primastimad}, we obtain the estimate
\Bsist
  && (\munp,\ynp-\mz)  
  = (\munp-\mean\munp , \ynp-\mz) 
  + (\mean\munp , \ynp-\mz)
  \non
  \\[1mm]
  && \leq \,\hat c\, \norma{A^r\munp} \, \norma{\ynp-\mz} 
  + (\mean\munp , -h\munp)
  \non
  \\[1mm]
  && \leq\,  c \, \norma{A^r\munp}
  - |\Omega| \, h \, (\mean\munp)^2
  \leq c \, \norma{A^r\munp} \,.
  \non
\Esist
\gianni{Therefore}, \eqref{perterzastima} becomes
\Beq
  \norma{\betal(\ynp)}_{\Luno}
  \leq c \, \Bigl(
    \tau \, \Norma\dhyn
    + \norma\unp
    + \norma{A^r\munp}
    + 1
  \Bigr).
  \label{stimabetaln}
\Eeq
Now, we square \eqref{stimabetaln}, multiply by~$h$
and sum up over~$n=0,\dots,k-1$ with $k\leq N$.
We deduce that
\begin{align*}
  &\somma n0{k-1} h \norma{\betal(\ynp)}_{\Luno}^2
  \\
  &\leq \,c \, \tau \, h \somma n0{k-1} \Norma\dhyn^2
  + c \, h \somma n0{k-1} \norma\unp^2
  + c \, h \somma n0{k-1} \norma{A^r \munp}^2 
  + c \,.
  \non
\end{align*}
Thanks to \eqref{primastimad}, the \rhs\ is bounded, and we conclude that
\Beq
  \norma{\betal(\overyh)}_{\L2\Luno} \leq c \,.
  \label{stimabetal}
\Eeq
At this point, we simply integrate \eqref{secondah} over~$\Omega$ 
and have \aet
\Bsist
  & |\Omega| \mean\overmuh
  & = \,\tau \iO \dt\yh
  + \Lpi' \iO (\overyh-\underyh)
  + \bigl( B^\sigma\overyh,B^\sigma(1) \bigr)
  \non
  \\
  && \quad {}
  + \iO \betal(\overyh)
  + \iO \pi(\overyh)
  - \iO \overuh \,.
  \non
\Esist
Thus, $\,\mean\overmuh\,$ is bounded in $L^2(0,T)$, thanks to \eqref{primastima}, \eqref{stimabetal} and~\eqref{hpu}.
This completes the proof of the desired estimate~\eqref{terzastima} as far as 
$\,\overmuh\,$ is concerned.
Since $A^r\mu^0=A^r0=0$ and $\,\overmuh-\undermuh\,$ is bounded in~$\L2H$ 
by virtue of \eqref{primastima}, 
the same estimate holds for~$\,\undermuh$.
\gianni{Hence, \eqref{terzastima} holds also in the case $\lambda_1=0$.}

\step
Limit

Collecting the estimates \accorpa{primastima}{terzastima},
and using standard weak and weak-star compactness results, we see that there are functions
$\yl$ and $\mul$ such that
\Bsist
  && \overyh \to \yl \,, \quad
  \underyh \to \yl\,,
  \aand
  \yh \to \yl
  \quad \hbox{weakly star in $\L\infty{\VB\sigma}$},
  \label{convyh}
  \\
  && \dt\yh \to \dt\yl
  \quad \hbox{weakly in $\L2{\VA{-r}}$},
  \label{convdtyh}
  \\[1mm]
  && \dt\yh \to \dt\yl
  \quad \hbox{weakly in $\L2H$\quad if $\tau>0$},
  \label{convdtyhbis}
  \\[1mm]
  && \overmuh \to \mul
  \quad \hbox{weakly in $\L2{\VA r}$},
  \label{convmuh}
\Esist
as $h\searrow0$ \,(more precisely, as $N\to\infty$), 
at least for some (not relabeled) subsequence, 
provided that $\lambda>0$ is small enough.
By letting $h$ tend to zero in~\eqref{cauchyh}, 
we see that $\yl$ satisfies~\eqref{cauchyl}.
Now, we prove that
\Beq
  \undermuh \to \mul
  \quad \hbox{weakly in $\L2{\VA r}$}.
  \label{convundermuh}
\Eeq
By \eqref{primastima} and \eqref{convmuh}, it suffices to check that
\Beq
  {}_{\L2{\VA{-r}}}\<v , \overmuh-\undermuh>_{\L2{\VA r}} \to 0
  \quad \hbox{as $h\searrow0$},
  \label{perconv}
\Eeq
for every $v$ belonging to a dense subspace $\calS$ of $\L2{\VA{-r}}$,
where we can take $\calS=\pier{C^1_c}(0,T;H)$ since $H$ is dense in~$\VA{-r}$
(see~\eqref{compembAneg}).
So, we fix $v\in \pier{C^1_c}(0,T;H)$ and choose $\delta>0$ such that 
$v(t)=0$ for $t\in[0,T]\setminus(\delta,T-\delta)$.
If $h\in(0,\delta/2)$, then we have
\begin{align}
  & |{}_{\L2{\VA{-r}}}\<v , \overmuh-\undermuh>_{\L2{\VA r}}|
  = \Bigl| \int_h^T \bigl((\overmuh-\undermuh)(t) , v(t) \bigr)\, dt \Bigr|
  \non
  \\
  & = \Bigl| \int_h^T \bigl(  \overmuh(t) - \overmuh(t-h)  , v(t) \bigr)\, dt \Bigr|
  = \Bigl| \int_h^T \bigl(\overmuh(t) , v(t) \bigr)\, dt
  - \int_0^{T-h} \bigl(\overmuh(t) , v(t+h) \bigr)\, dt \Bigr|
  \non
  \\
  & = \Bigl| \int_h^{T-h} \bigl(\overmuh(t) , v(t)-v(t+h)\bigr) \, dt \Bigr|
  \leq \norma\overmuh_{\L2H}\, \norma{v'}_{\L\infty H} \, h^{1/2}\,,
  \non
\end{align}
and \eqref{perconv} follows.
Therefore, \eqref{convundermuh} is proved and the pair $(\yl,\mul)$ solves~\eqref{primal}.
In order to deal with the nonlinear terms of~\eqref{secondah}, 
we owe to the compact embedding $\VB\sigma\subset H$ (see~\eqref{compembB})
and to well-known strong compactness results (see, e.g., \cite[Sect.~8, Cor.~4]{Simon}).
From \accorpa{convyh}{convdtyh} we deduce that
\Beq
  \yh \to \yl
  \quad \hbox{strongly in $\L\infty H$}.
  \label{strongyh}
\Eeq
\gianni{This and \eqref{primastima} (see the last term on the \lhs) imply that}
\Beq
  \overyh \to \yl
  \quad \hbox{strongly in $\L2H$}.
  \label{strongoveryh}
\Eeq
By \Lip\ continuity, we infer that also
\Beq
  (\betal+\pi)(\overyh) \to (\betal+\pi)(\yl)
  \quad \hbox{strongly in $\L2H$}.
  \non
\Eeq
Moreover, as we can assume that $\overyh$ converges to $\yl$ \aeQ\
and $\Betal$ grows at most quadratically, we can also apply \pier{\eqref{daprimastima} and Fatou's lemma} to deduce that
\Beq
  \iO\Betal(\yl(t))
  \pier{{}\leq  \liminf_{h\searrow 0}}\iO\Betal(\overyh(t))
  \leq c 
  \quad \aat ,
\Eeq
whence
\Beq
  \norma{\Betal(\yl)}_{\L\infty\Luno}
  \leq c \,.
  \label{stimaBetalyl}
\Eeq
Therefore, we can pass to the limit in the time-integrated version of \eqref{secondah}
(written with time-dependent test functions)
and deduce that the pair $(\yl,\mul)$ also solves \eqref{intsecondal},
which is equivalent to~\eqref{secondal}.
This concludes the proof of Theorem~\ref{Wellposednessl}.


\section{Existence}
\label{EXISTENCE}
\setcounter{equation}{0}

This section is devoted to the proof of the existence part
of Theorem~\ref{Wellposedness}.
Just by the semicontinuity of the norms,
all of the uniform estimates we have established 
for the interpolants of the discrete solution
hold true for the (unique) solution to the approximating problem.
Therefore, from \accorpa{primastima}{secondastima}, 
\eqref{terzastima} and~\eqref{stimaBetalyl}, 
we deduce~that 
\Bsist
  && \norma\yl_{\H1{\VA{-r}}\cap\L\infty{\VB\sigma}}
  + \norma\mul_{\L2{\VA r}}
  \non
  \\[1mm]
  && 
  + \tau^{1/2} \norma{\dt\yl}_{\L2H}
  + \norma{\Betal(\yl)}_{\L\infty\Luno}
  \leq c 
  \label{stimasoluzl}
\Esist
for $\lambda>0$ small enough.
We infer that 
there exist a strictly decreasing sequence $\lambda_n\searrow0$ 
and a pair $(y,\mu)$ satisfying, as~$n\nearrow\infty$,
\Bsist
  && \yln \to y
  \quad \hbox{weakly star in $\H1{\VA{-r}}\cap\L\infty{\VB\sigma}$}\,,
  \label{convyl}
  \\
  && \muln \to \mu
  \quad \hbox{weakly in $\L2{\VA r}$}\,,
  \label{convmul}
  \\
  && \dt\yln \to \dt y
  \quad \hbox{weakly in $\L2H$\quad if $\tau>0$} \,.
  \label{convdtyl}
\Esist
Then, it is immediately seen that $(y,\mu)$ solves \eqref{prima}
and that $y$ satisfies the initial condition~\eqref{cauchy}.
Now, we prove that the variational inequality \eqref{seconda} holds true as well.
To this end, we first owe to the compact embedding $\VB\sigma\subset H$ (see~\eqref{compembB})
and, e.g., to \cite[Sect.~8, Cor.~4]{Simon}),
and deduce that we also have\pier{, at least for another subsequence, that} 
\Beq
  {\yln} \to y
  \quad \hbox{strongly in $\L\infty H$}
  \aand
  \aeQ .
  \label{strongyl}
\Eeq
This implies that $\pi({\yln})$ converges 
to $\pi(y)$ in the same space, by \Lip\ continuity.
Now, we use \eqref{strongyl} once more to show that
\Beq
  \intQ \Beta(y)
  \leq \liminf_{n\to\infty} \intQ \Betaln(\yln)
  < +\infty .
  \label{liminfBeta}
\Eeq
We notice that the \rhs\ of \eqref{liminfBeta} actually is finite 
thanks to the bound for $\Betal(\yl)$ given by~\eqref{stimasoluzl}.
In particular, the requirement $\Beta(y)\in\LQ1$ (see~\eqref{regBetay}) is fulfilled
once the first inequality of \eqref{liminfBeta} is established.
In order to prove~it,
we take arbitrary indices $m$ and $n$ with $n>m$.
Then\, $\lambda_n<\lambda_m$, and we can apply~\eqref{propBetal}.
We deduce that
\Beq
  \Betalm(\yln)
  \leq \Betaln(\yln)
  \quad \hbox{\aeQ, for every $n>m$},
  \non
\Eeq
whence also (since $\Betalm$ is continuous)
\Beq
  \Betalm(y)
  = \lim_{n\to\infty} \Betalm(\yln)
  = \liminf_{n\to\infty} \Betalm(\yln)
  \leq \liminf_{n\to\infty} \Betaln(\yln)
  \quad \aeQ.
  \non
\Eeq
On the other hand, we have, by {virtue of the second
property stated in}~\eqref{propBetal},
\Beq
  \Beta(y) = \lim_{m\to\infty} \Betalm(y)
  \quad \aeQ.
\Eeq
Thus,
\Beq
  \Beta(y)
  \leq \liminf_{n\to\infty} \Betaln(\yln)
  \quad \aeQ ,
\Eeq
and \eqref{liminfBeta} follows from Fatou's lemma.
Next, we have that
\Bsist
  && \ioT \bigl( B^\sigma y(t) , B^\sigma (y(t) - v(t)) \bigr) \, dt
  \non
  \\
  && \leq \liminf_{n\to\infty} \ioT \bigl( B^\sigma \yln(t) , B^\sigma \yln(t) \bigr) \, dt
  - \lim_{n\to\infty} \ioT \bigl( B^\sigma \yln(t) , B^\sigma v(t) \bigr) \, dt
  \non
  \\
  && = \liminf_{n\to\infty} \ioT \bigl( B^\sigma \yln(t) , B^\sigma (\yln(t) - v(t)) \bigr) \, dt
  \non
\Esist
for every $v\in\L2{\VB\sigma}$,
since $B^\sigma\yl$ converges to $B^\sigma y$ weakly in $\L2H$ by \eqref{convyl}.
At this point, we can let $n$ tend to infinity in \eqref{intsecondal}
 written with $\lambda=\lambda_n$.
By also accounting for \eqref{convmul}, \eqref{strongyl} and~\eqref{propBetal},
we see that, for every $v\in\L2{\VB\sigma}$, we~have
\begin{align}
  & \intQ \Beta(y)
  + \ioT \bigl( B^\sigma y(t) , B^\sigma (y(t) - v(t)) \bigr) \, dt
  \non
  \\
  & \leq \,\liminf_{n\to\infty} \intQ \Betaln(\yln)
  + \liminf_{n\to\infty} \ioT \bigl( B^\sigma \yln(t) , B^\sigma (\yln(t) - v(t)) \bigr) \, dt
  \non
  \\
  & \leq \,\liminf_{n\to\infty} \Bigl(
    \intQ \Betaln(\yln)
    + \ioT \bigl( B^\sigma \yln(t) , B^\sigma (\yln(t) - v(t)) \bigr) \, dt
  \Bigr)
  \non
  \\
  & \leq \,\lim_{n\to\infty} \,\Bigl( \ioT \bigl(
\pier{{}- \tau \partial_t y^{\lambda_n}(t)} - \pi(\yln(t)) + u(t) \juerg{\,+\,} \muln(t) , \yln(t) - v(t)
  \bigr)  \, \pier{dt}
  + \intQ \Betaln(v)\Bigr)
  \non
  \\
  & = \,\ioT \bigl(\pier{{}-\tau \partial_t y(t)} 
     -\pi(y(t)) + u(t) {\,+\,} \mu(t) , y(t)-v(t) \bigr) \, dt
    + \intQ \Beta(v).
  \non
\end{align}
Thus, \eqref{intseconda} holds true.
Since \eqref{intseconda} is equivalent to~\eqref{seconda}, 
the proof of Theorem~\ref{Wellposedness} is complete.


\section{Regularity}
\label{REGULARITY}
\setcounter{equation}{0}

This section is devoted to the proof of Theorem~\ref{Regularity}.
Coming back to the proofs of Theorems~\ref{Wellposedness} and~\ref{Wellposednessl},
we see that it is sufficient to establish some estimates 
on the solution to the discrete problem in one of the forms
\accorpa{primad}{secondad} and \accorpa{primah}{secondah},
uniformly with respect to both $h$ and~$\lambda$.
Of course, we can account for the estimates proved in Section~\ref{APPROX}.

\step
First regularity estimate

We prove the uniform estimate
\Beq
  \norma{A^r\overmuh}_{\L\infty H}
  + \norma{B^\sigma\dt\yh}_{\L2H}
  + \tau^{1/2} \norma{\dt\yh}_{\L\infty H}
  \leq c\,,
  \label{stimaregd}
\Eeq
with a constant $c$ that does not depend on $h$, $\lambda$ and~$\tau$
(like the constant $K_3$ in the statement of the theorem).
We test \eqref{primad} by $\munp-\mun$.
On account of~\eqref{propA}, we obtain
\Beq
  \Bigl( \dhyn , \munp-\mun \Bigr)
  + \norma{\munp-\mun}^2
  + \bigl( A^r \munp , A^r(\munp-\mun) \bigr) = 0.
  \qquad
  \label{testprimad}
\Eeq
Now, we perform a discrete differentiation on~\eqref{secondad}.
Precisely, we write it for both $(\yn,\mun)$ and $(\ynm,\munm)$, 
take the difference, divide by~$h$ and rearrange.
We have for $1\leq n<N$
\Bsist
  && \frac 1h \Bigl( \tau \, \dhyn - \tau \, \dhynm \Bigr)
  + \Lpi' \Bigl( \dhyn - \dhynm \Bigr)
  \non
  \\
  && \quad {}
  + B^{2\sigma} \, \dhyn
  + \frac 1h \bigl( \betal(\ynp) - \betal(\yn) \bigr)
  \non
  \\
  && = \dhmun
  + \dhun
  - \frac 1h \bigl( \pi(\ynp) - \pi(\yn) \bigr)
  \non
\Esist
and test this equality by $\ynp-\yn$.
\pier{On account of} \eqref{propB}, we obtain
\Bsist
  && \tau \Bigl( \dhyn , \dhyn-\dhynm \Bigr)
  \non
  \\
  && \quad {}
  + \Lpi' h \Bigl( \dhyn , \dhyn-\dhynm \Bigr)
  \non
  \\
  && \quad {}
  + h \, \Norma{B^\sigma \, \dhyn}^2
  + \frac 1h \bigl( \betal(\ynp) - \betal(\yn) , \ynp-\yn \bigr)
  \non
  \\
  && = \Bigl( \dhmun , \ynp-\yn \Bigr)
  \non
  \\
  && \quad {}
  + \Bigl( \dhun , \ynp-\yn \Bigr)
  - \frac 1h \bigl( \pi(\ynp) - \pi(\yn) , \ynp-\yn \bigr).
  \label{testdsecondad}
\Esist
Now, we add this to \eqref{testprimad} and notice that two terms cancel each other
and that the term involving $\betal$ is nonnegative by monotonicity.
Thus, \pier{thanks} to the identity~\eqref{elementare},
applying the Schwarz and Young inequalities to the remaining terms on the \rhs,
and accounting for the \Lip\ continuity of~$\pi$,
we deduce that
\Bsist
  && \norma{\munp-\mun}^2
  + \frac 12 \norma{A^r \munp}^2
  + \frac 12 \, \norma{A^r(\munp-\mun)}^2
  - \frac 12 \, \norma{A^r \mun}^2
  \non
  \\
  && \quad {}
  + \frac \tau 2 \, \Bigl( \Norma\dhyn^2 - \Norma\dhynm^2 \Bigr)
  + \frac \tau 2 \, \Norma{\dhyn-\dhynm}^2
  \non
  \\
  && \quad {}
  + \Lpi' \, \frac h2 \, \Bigl( \Norma\dhyn^2 - \Norma\dhynm^2 \Bigr)
  \non
  \\
  && \quad {}
  + \Lpi' \, \frac h2 \, \Norma{\dhyn-\dhynm}^2
  + h \, \Norma{B^\sigma \, \dhyn}^2
  \non
  \\
  && \leq \frac h2 \, \Norma\dhun^2
  + \frac h2 \, \Norma\dhyn^2
  + \frac {\Lpi h} 2 \, \Norma\dhyn^2 \,.
  \non
\Esist
Summing up for $n=1,\dots,k-1$ with $k\leq N$, 
and omitting a number of nonnegative terms on the \lhs, we infer that
\Bsist
  && \frac 12 \norma{A^r \muk}^2
  + \frac \tau 2 \,  \Norma{\frac {y^k-y^{k-1}} h}^2
  + \somma n1{k-1} h \, \Norma{B^\sigma \, \dhyn}^2
  \non
  \\
  && \leq \frac 12 \, \norma{A^r \muu}^2
  + \frac \tau 2 \, \Norma{\frac {\yu-\yz} h}^2
  + \Lpi' \, \frac h2 \, \Norma{\frac {\yu-\yz} h}^2
  \non
  \\
  && \quad {}
  + \somma n1{k-1} h \Norma\dhun^2
  + \frac {\Lpi'} 2 \somma n1{k-1} h \Norma\dhyn^2 \,.
  \non
\Esist
At this point, we use~\eqref{interpH1Z}, the compactness inequality \eqref{compineq}
and the estimate \eqref{secondastima}, to \pier{control the last two terms on the \rhs:}
\Bsist
  && \somma n1{k-1} h \Norma\dhun^2
  + \frac {\Lpi'} 2 \somma n1{k-1} h \Norma\dhyn^2
  \non
  \\
  && \leq \norma{\dt u}_{\L2H}^2
  + \frac 12 \somma n1{k-1} h \, \Norma{B^\sigma \, \dhyn}^2
  + c \somma n1{k-1} h \, \Norma{\dhyn}_{A,-r}^2
  \non
  \\
  && = c
  + \frac 12 \somma n1{k-1} h \, \Norma{B^\sigma \, \dhyn}^2
  + c \, \norma{\dt\yh}_{\L2{A,-r}}^2
  \leq \frac 12 \somma n1{k-1} h \, \Norma{B^\sigma \, \dhyn}^2
  + c \,.
  \non
\Esist
Therefore, on account of Proposition~\ref{Propinterp}, the above inequality becomes
\Bsist
  && \norma{A^r \overmuh}_{\L\infty H}^2
  + \tau \, \Norma{\dt\yh}_{\L\infty H}^2
  + \norma{B^\sigma \dt\yh}_{\L2H}^2
  \non
  \\
  && \leq c \Bigl( 
    \norma{A^r \muu}^2
    + \tau \, \Big\|\frac {\yu-\yz} h\Big\|^2
    + h \, \Big\|\frac {\yu-\yz} h\Big\|^2
    + 1
  \Bigr)\,,
  \label{perstimareg}
\Esist
and \eqref{stimaregd} will follow whenever we estimate the \rhs\ of~\eqref{perstimareg}.
To this end,
we write \eqref{primad} and \eqref{secondad} with $n=0$.
We also rearrange the latter, recall that $y^0=\yz$ and $\mu^0=0$, 
and set for convenience $\calA_\lambda:=\Lpi' I + \betal + \pi$.
We~have
\Bsist
  && \frac {\yu-\yz} h + \muu + A^{2r} \muu
  = 0
  \label{primadz}
  \\
  && \tau \, \frac {\yu-\yz} h
  + \calA_\lambda (\yu)
  - \calA_\lambda (\yz)
  + B^{2\sigma} (\yu-\yz)
  \non
  \\
  && = \muu
  + \bigl(
    u^1 
    - B^{2\sigma} \yz
    - \betal(\yz)
    - \pi(\yz)  
  \bigr).
  \label{secondadz}
\Esist
Now, we test \eqref{primadz} by $\muu$ 
and \eqref{secondadz} by $(\yu-\yz)/h=-(\muu+A^{2r}\muu)$,
by choosing the first or second expression according to our convenience.
\pier{In view of}~\accorpa{propA}{propB},
and noting that $(\calA_\lambda(y^1)-\calA_\lambda(\yz),y^1-\yz)\geq0$
since $\betal$ is monotone and $\Lpi$ is the \Lip\ constant of~$\pi$, 
we obtain
\Beq
  \Bigl( \frac {\yu-\yz} h , \muu \Bigr)
  + \norma\muu^2
  + \norma{A^r \muu}^2
  = 0
  \label{testprimadz}
\Eeq
\pier{first, and then}   
  \begin{align}
  & \tau \Norma{\frac{\yu-\yz}h}^2
  + h \, \norma{B^\sigma(\yu-\yz \pier{)}}^2
  \non
  \\
  & \quad {} \leq \Bigl( \muu , \frac {\yu-\yz} h \Bigr)
  + \Bigl(
    u^1 
    - B^{2\sigma} \yz
    - \betal(\yz)
    - \pi(\yz) ,
    \frac {\yu-\yz} h
  \Bigr)
  \label{testsecondadzy}
\end{align} 
\pier{or, alternatively,}
\begin{align}
  & \tau \Norma{\frac{\yu-\yz}h}^2
  + h \, \norma{B^\sigma(\yu-\yz)}^2
  \non
  \\
  & \quad {} \leq \Bigl( \muu , \frac {\yu-\yz} h \Bigr)
  - \bigl(
    u^1 
    - B^{2\sigma} \yz
    - \betal(\yz)
    - \pi(\yz) ,
    \muu + A^{2r}\muu
  \bigr)\,.
  \label{testsecondadzmu}
\end{align}
Now, we distinguish the two cases of the statement of Theorem~\ref{Regularity}.
We first assume $\tau>0$ and \eqref{hpyzreg}.
Then, we add \eqref{testprimadz} and \eqref{testsecondadzy},
by noticing that two terms cancel each other.
Moreover, we omit a nonnegative term on the \lhs\ and use the Schwarz and Young inequalities on the \rhs.
We then have that
\Beq
  \norma\muu^2
  + \norma{A^r \muu}^2
  + \tau \Norma{\frac{\yu-\yz}h}^2
  \leq \frac \tau 2 \, \Norma{\frac {\yu-\yz} h}^2
  + \frac 1 {2\tau} \norma{u^1-B^{2\sigma}\yz-\betal(\yz)-\pi(\yz)}^2 .
  \non
\Eeq
By accounting for \eqref{hpu}, which implies that\, $\norma{u^1}=\norma{u(h)}\leq c$,
\eqref{hpyzreg} and~\eqref{propbetal}, we see that the last norm is bounded
uniformly with respect to~$\lambda$.
Therefore, the \rhs\ of \eqref{perstimareg} is bounded, too.
Now, we assume $\tau=0$ and \accorpa{hpdatireg}{defmul}.
Then, we add \eqref{testprimadz} and \eqref{testsecondadzmu} and similarly have that
\Bsist
  && \norma\muu^2
  + \norma{A^r \muu}^2
  \leq - \bigl(
    u^1 
    - B^{2\sigma} \yz
    - \betal(\yz)
    - \pi(\yz) ,
    \muu + A^{2r}\muu
  \bigr)
  \non
  \\
  && \leq \frac 12 \, \norma\muu^2
  + \frac 12 \, \norma{u^1 - B^{2\sigma} \yz - \betal(\yz) - \pi(\yz)}^2
  \non
  \\
  && \quad {}
  + \frac 12 \, \norma{A^r \muu}^2
  + \frac 12 \, \norma{A^r(u^1 - B^{2\sigma} \yz - \betal(\yz) - \pi(\yz))}^2 .
  \non
\Esist
Hence, the sought bound is ensured by \accorpa{hpdatireg}{defmul},
since $u^1=u(h)$.
Therefore, \eqref{stimaregd} is established in any case.

\step
Consequence

By applying the compactness inequality~\eqref{compineq}, we obtain
\Beq
  \norma{\dt\yh(t)}^2
  \leq \norma{B^\sigma \dt\yh(t)}^2
  + c \, \norma{\dt\yh(t)}_{A,-r}^2
  \quad \aat .
  \non
\Eeq
On the other hand, we have that \,$\norma{\dt\yh}_{\L2{\VA{-r}}}\leq c$,
 by virtue of~\eqref{secondastima}.
Therefore, we deduce from \eqref{stimaregd} that
\Bsist
  && \norma{\dt\yh}_{\L2H}^2
  \leq \norma{B^\sigma\dt\yh}_{\L2H}^2
  + c \, \norma{\dt\yh}_{\L2{\VA{-r}}}^2
  \leq c \,,
  \non
  \\
  \noalign{\noindent as well as}
  && \norma{\dt\yh}_{\L2{\VB\sigma}}^2
  = \norma{\dt\yh}_{\L2H}^2
  + \norma{B^\sigma\dt\yh}_{\L2H}^2
  \leq c \,.
  \non
\Esist
This implies that
\Beq
  \dt y \in \L2{\VB\sigma}
  \aand
  \norma{\dt y}_{\L2{\VB\sigma}}
  \leq c\,,
  \non
\Eeq
which is a part of \eqref{regsoluzbis} and~\eqref{stimareg}.

\step
Second regularity estimate

We now prove the inequalities
\Beq
  \norma\overmuh_{\L\infty H}
  \leq c
  \aand
  \norma\undermuh_{\L\infty H}
  \leq c\,,
  \label{stimahLinfty}
\Eeq
the latter being a consequence of the former since $\mu^0=0$.
If $\lambda_1>0$, then\, $\norma v\leq \norma{A^rv}$ for every $v\in\VA r$,
so that \eqref{stimaregd} also implies that
\Beq
  \norma\overmuh_{\L\infty H}
  \leq c \, \norma{A^r\overmuh}_{\L\infty H}\,
  \leq c\, ,
  \non
\Eeq
and the {first claim of \eqref{stimahLinfty} is proved for the case $\lambda_1>0$.
In the case \,}$\lambda_1=0$, we only have
(see \eqref{defnormaAr} and Remark~\ref{RemnormaVAr})
\Beq
  \norma{\overmuh-\mean\overmuh}_{\L\infty H}
  \leq c \, \norma{A^r\overmuh}_{\L\infty H}
  \leq c \,.
  \non
\Eeq
Thus, in order to achieve \eqref{stimahLinfty}, we have to estimate the mean value.
To this end, we recall \eqref{stimabetaln}, which can be written in the form
\Beq
  \norma{\betal(\overyh(t))}_{\Luno}
  \leq c \bigl(
    \tau \, \norma{\dt\yh(t)}
    + \norma{\overuh(t)}
    + \norma{A^r \overmuh(t)}
    + 1
  \bigr)
  \leq c
  \quad \aat .
  \non
\Eeq
From \eqref{stimaregd} and~\eqref{hpu}, we deduce that
\Beq
  \norma{\betal(\overyh)}_{\L\infty\Luno}
  \leq c \,.
  \non
\Eeq
At this point, we simply integrate \eqref{secondah} over~$\Omega$ 
to obtain, almost everywhere in $(0,T)$,
\Bsist
  & |\Omega| \mean\overmuh
  & = {\tau} \iO \dt\yh
  + \Lpi' \iO (\overyh-\underyh)
  + \bigl( B^\sigma\overyh,B^\sigma(1) \bigr)
  \non
  \\
  && \quad {}
  + \iO \betal(\overyh)
  + \iO \pi(\overyh)
  - \iO \overuh \,.
  \non
\Esist
Thus, $\mean\overmuh$ is bounded in $L^\infty(0,T)$ thanks to \eqref{primastima} 
and~\eqref{hpu}.
This concludes the proof of~\eqref{stimahLinfty}.

\step
Conclusion

From \eqref{stimaregd} and \eqref{stimahLinfty}, we infer that
\Beq
  \mu \in \L\infty{\VA r}
  \aand
  \norma\mu_{\L\infty{\VA r}} \leq c\,,
  \non
\Eeq
which is another claim of \eqref{regsoluzbis} and~\eqref{stimareg}.
Moreover, by recalling \eqref{primah} and Proposition~\ref{ExtensionA},
we deduce that
\Bsist
  && \norma{\dt\yh}_{\L\infty{\VA{-r}}}
  \leq \norma{A^{2r}\muh}_{\L\infty{\VA{-r}}}
  + c \, \norma{\undermuh-\overmuh}_{\L\infty H}
  \non
  \\
  && \leq \norma{A^r\muh}_{\L\infty H} + c 
  \,\leq \,c \,,
  \non
\Esist
which yields that
\Beq
  \dt y \in \L\infty{\VA{-r}}
  \aand
  \norma{\dt y}_{\L\infty{\VA{-r}}} \leq c \,.
  \non
\Eeq
Now, we assume that $\tau>0$, in addition.
Then \eqref{secondah}, \eqref{stimareg} and \eqref{stimahLinfty} give that 
\Beq
  \tau^{1/2} \norma{A^{2r}\overmuh}_{\L\infty H}
  \leq \tau^{1/2} \norma{\dt\yh}_{\L\infty H}
  + \tau^{1/2} \norma{\undermuh-\overmuh}_{\L\infty H}
  \leq c \,,
  \non
\Eeq
whence
\Beq
  \dt y \in \L\infty H, \enskip
  \mu \in \L\infty{\VA{2r}}
  \aand
  \norma{\tau^{1/2}\dt y}_{\L\infty H}
  + \norma{\tau^{1/2}\mu}_{\L\infty{\VA{2r}}}
  \leq c\,.
  \non
\Eeq
This concludes the proof of Theorem~\ref{Regularity}.


\section*{Acknowledgments}
The research of PC is supported by the Italian Ministry of Education, 
University and Research~(MIUR): Dipartimenti di Eccellenza Program (2018--2022) 
-- Dept.~of Mathematics ``F.~Casorati'', University of Pavia. In addition, PC and GG gratefully acknowledge some financial support \pier{from 
the MIUR-PRIN Grant 2015PA5MP7 ``Calculus of Variations'',}
the GNAMPA (Gruppo Nazionale per l'Analisi Matematica, 
la Probabilit\`a e le loro Applicazioni) of INdAM (Isti\-tuto 
Nazionale di Alta Matematica) and the IMATI -- C.N.R.~Pavia.  


\vspace{3truemm}

\Begin{thebibliography}{10}

\pier{%
\bibitem{ABG}
H. Abels, S. Bosia, M. Grasselli, Cahn--Hilliard equation with nonlocal singular free energies, Ann. Mat. Pura Appl. (4)\/
{\bf 194} (2015), 1071-1106.
\bibitem{AD}
N. Abatangelo, L. Dupaigne,
Nonhomogeneous boundary conditions for the spectral fractional Laplacian,
Ann. Inst. H. Poincar\'e Anal. Non Lin\'eaire {\bf 34} (2017), 439-467.
\bibitem{AM}
M. Ainsworth, Z. Mao, Analysis and approximation of a fractional
Cahn--Hilliard equation, SIAM J. Numer. Anal. {\bf 55} (2017), 1689-1718.
\bibitem{AkSS1}
G. Akagi, G. Schimperna, A. Segatti,
Fractional Cahn--Hilliard, Allen--Cahn and porous medium equations,
J. Differential Equations {\bf 261} (2016), 2935-2985.
\bibitem{AkSS2}
G. Akagi, G. Schimperna, A. Segatti,
Convergence of solutions for the fractional Cahn--Hilliard system,
J. Funct. Anal. {\bf 276} (2019), 2663-2715.}
\bibitem{Barbu}
V. Barbu,
``Nonlinear Differential Equations of Monotone Type in Banach Spaces'',
Springer,
London, New York, 2010.

\newpier{%
\bibitem{BSX}
Z. Binlin, M. Squassina, Z. Xia,
Fractional {NLS} equations with magnetic field, 
critical frequency and critical growth,
Manuscripta Math. {\bf 155} (2018), 115-140.}

\pier{%
\bibitem{BFV}
M. Bonforte, A. Figalli, J.L. V\'azquez, 
Sharp global estimates for local and nonlocal porous medium-type equations in bounded 
domains, Anal. PDE {\bf 11} (2018), 945-982. 
\bibitem{BSV}
M. Bonforte, Y. Sire, J.L. V\'azquez,  
Existence, uniqueness and asymptotic behaviour for fractional porous medium on bounded domains, 
Discrete Contin. Dyn. Syst. {\bf 35} (2015), 5725-5767.
\bibitem{BV}
M. Bonforte, J.L. V\'azquez,  
A priori estimates for fractional nonlinear 
degenerate diffusion equations on bounded domains, 
Arch. Ration. Mech. Anal. {\bf 218} (2015), 317-362.}

\newpier{%
\bibitem{BSY}
L. Brasco, M. Squassina, Y. Yang,
Global compactness results for nonlocal problems,
Discrete Contin. Dyn. Syst. Ser. S\/ {\bf 11} (2018), 391-424.}

\bibitem{Brezis}
H. Brezis,
``Op\'erateurs maximaux monotones et semi-groupes de contractions
dans les espaces de Hilbert'',
North-Holland Math. Stud.
{\bf 5},
North-Holland,
Amsterdam,
1973.

\pier{%
\bibitem{CT}
X. \newpier{Cabr\'e},  J. Tan, 
Positive solutions of nonlinear problems involving the square root of the 
{L}aplacian, Adv. Math. {\bf 224} (2010), 2052-2093.
\bibitem{CahH} 
J.W. Cahn, J.E. Hilliard, 
Free energy of a nonuniform system I. Interfacial free energy, 
J. Chem. Phys. {\bf 2} (1958), 258-267.
\bibitem{CS}
L.A. Caffarelli, P.R. Stinga,
Fractional elliptic equations, {C}accioppoli estimates and regularity,
Ann. Inst. H. Poincar\'e Anal. Non Lin\'eaire {\bf 33} (2016), 767-807.
\bibitem{CF2} 
P.\ {C}olli, T.\ {F}ukao, 
Equation and dynamic boundary condition of 
Cahn--Hilliard type with singular potentials, 
Nonlinear Anal. {\bf 127} (2015), 413-433.
\bibitem{CGS6}
P. Colli, G. Gilardi, J. Sprekels,
On an application of Tikhonov's fixed point theorem
to a nonlocal Cahn--Hilliard type system modeling phase separation,
J. Differential Equations {\bf 260} (2016), 7940-7964.
\bibitem{CGS13}
P. Colli, G. Gilardi, J. Sprekels,
On a Cahn--Hilliard system with convection and
dynamic boundary conditions, 
Ann. Mat. Pura Appl. (4)\/ {\bf 197} (2018), 1445-1475.
\bibitem{DPV}
E. Di Nezza, G. Palatucci, E. Valdinoci,
Hitchhiker's guide to the fractional Sobolev spaces,
Bull. Sci. Math. {\bf 136} (2012), 521-573.
\bibitem{DG}
M. D'Elia, M. Gunzburger, 
The fractional Laplacian operator on bounded domains as a special case of the
nonlocal diffusion operator,
Comput. Math. Appl. {\bf 66} (2013), 1245-1260.
\bibitem{EZ} 
C.M. Elliott, S. Zheng, 
On the Cahn--Hilliard equation, 
Arch. Rational Mech. Anal. {\bf 96} (1986), 339-357.
\bibitem{FG} 
E. Fried, M.E. Gurtin, 
Continuum theory of thermally induced phase transitions based on an order 
parameter,  Phys. D {\bf 68} (1993), 326-343.
\bibitem{GalDCDS}
C.G. Gal, 
On the strong-to-strong interaction case 
for doubly nonlocal Cahn--Hilliard equations,
Discrete Contin. Dyn. Syst.  {\bf 37} (2017),  131-167.
\bibitem{GalEJAM}
C.G. Gal, 
Non-local Cahn--Hilliard equations with fractional dynamic boundary conditions, European J. Appl. Math. {\bf 28} (2017), 736-788. 
\bibitem{GalAIHP}
C.G. Gal, 
Doubly nonlocal Cahn--Hilliard equations, 
Ann. Inst. H. Poincar\'e Anal. Non Lin\'eaire {\bf 35} (2018), 357-392.
\bibitem{GGG}
C.G. Gal,  A. Giorgini, M. Grasselli,  
The nonlocal Cahn--Hilliard equation with singular potential: well-posedness, regularity and strict separation property,
J. Differential Equations  {\bf 263} (2017), 5253-5297.}
\bibitem{GiMiSchi} 
G. Gilardi, A. Miranville, G. Schimperna,
On the Cahn--Hilliard equation with irregular potentials and dynamic boundary conditions,
Commun. Pure Appl. Anal. {\bf 8} (2009), 881-912.
\pier{%
\bibitem{Gru1} 
G. Grubb,  
Regularity of spectral fractional Dirichlet and Neumann problems,
Math. Nachr.  {\bf 289} (2016), 831-844.
\bibitem{Gru2} 
G. Grubb,  
Regularity in $L_p$ Sobolev spaces of solutions to fractional heat equations,
J. Funct. Anal.  {\bf 274}  (2018), 2634-2660.
\bibitem{Gu} 
M. Gurtin, 
Generalized Ginzburg-Landau and Cahn-Hilliard equations based on a microforce balance, Phys.~D {\bf 92} (1996), 178-192.}

\bibitem{Jerome}
J.W. Jerome, ``Approximation of Nonlinear Evolution Systems'', 
Math. Sci. Engrg.  {\bf 164}, Academic Press, Orlando,
1983.

\pier{%
\bibitem{Kwa}
M. Kwa\'snicki, Ten equivalent definitions of the fractional Laplace operator,
Fract. Calc. Appl. Anal. {\bf 20} (2017), 7-51.}

\bibitem{MiZe}
A. Miranville, S. Zelik, 
Robust exponential attractors for Cahn--Hilliard type equations with singular potentials, 
Math. Methods Appl. Sci. {\bf 27} (2004), 545-582.

\pier{%
\bibitem{MN}
R. Musina, A.I. Nazarov, 
Variational inequalities for the spectral fractional Laplacian,
Comput. Math. Math. Phys.  {\bf 57} (2017), 373-386.
\bibitem{NovCoh}
A.~Novick-Cohen, 
On the viscous {C}ahn--{H}illiard equation, in
``Material instabilities in continuum mechanics'' ({E}dinburgh, 1985--1986),
Oxford Sci. Publ., Oxford Univ. Press, New York, 1988, pp.~329-342.
\bibitem{RS1}
X. Ros-Oton, J. Serra, 
The Dirichlet problem for the fractional Laplacian: regularity up to
the boundary, J. Math. Pures Appl. (9) {\bf 101} (2014), 275-302.
\bibitem{RS2}
X. Ros-Oton, J. Serra, 
The extremal solution for the fractional Laplacian, 
Calc. Var. Partial Differential Equations {\bf 50} (2014), 723-750.
\bibitem{RH}
P. Rybka, K.-H. Hoffmann, 
Convergence of solutions to Cahn--Hilliard equation, 
Comm. Partial Differential Equations {\bf 24} (1999), 1055-1077.
\bibitem{SV0}
R. Servadei, E. Valdinoci, 
Variational methods for non-local operators of elliptic type, Discrete
Contin. Dyn. Syst. {\bf 33} (2013), 2105-2137.
\bibitem{SV1}		
R. Servadei, E. Valdinoci,  
On the spectrum of two different fractional operators,
Proc. Roy. Soc. Edinburgh Sect. A  {\bf 144}  (2014), 831-855.
\bibitem{SV2}		
R. Servadei, E. Valdinoci,  
Weak and viscosity solutions of the fractional Laplace equation,
Publ. Mat.  {\bf 58}  (2014), 133-154.
\bibitem{SV3}		
R. Servadei, E. Valdinoci, 
The Brezis-Nirenberg result for the fractional Laplacian,
Trans. Amer. Math. Soc.  {\bf 367} (2015), 67-102.}

\bibitem{Simon}
J. Simon,
Compact sets in the space $L^p(0,T; B)$,
Ann. Mat. Pura Appl.~(4)\/ 
{\bf 146}, (1987), 65-96.

\End{thebibliography}

\End{document}


  \\
  \separa
  && \iO \nabla w \cdot \nabla\xi
  + \iG \nablaG\wG \cdot \nablaG\xiG
  = \frac 12 \, \frac d{dt} \, \normaHH{(w,\wG)}^2
  \non
  \\
  && \quad \hbox{if $(w,\wG)\in\L2\calV$, \ $\dt(w,\wG) \in \L2\calVsz$}
  \non
  \\   
  && \aand
  \hbox{$\Xi = \calN(\dt(w,\wG))$} \,.
  && \< \gstar , \calN\gstar >_{\calV}
  = \norma\gstar_*^2
  \qquad \hbox{if $\gstar \in \calVsz$},
  \\
  \separa
  && \iO \nabla w \cdot \nabla\xi
  + \iG \nablaG\wG \cdot \nablaG\xiG
  = \normaHH{(w,\wG)}^2
  \non
  \\
  && \quad \hbox{if $(w,\wG) \in \calHz$ and $(\xi,\xiG) = \calN(w,\wG)$} \,.
  \\
  \separa